\nonstopmode
\documentclass[]{amsart}
\usepackage{amssymb}
\usepackage{amscd}
\usepackage{graphicx}
\usepackage{psfrag}
\usepackage{epsfig}
\usepackage[all]{xy}
\newdir{ >}{{}*!/-10pt/@{>}}
\newdir^{ (}{{}*!/-5pt/@^{(}}
\newdir_{ (}{{}*!/-5pt/@_{(}}
\def\cgaps#1{}
\def\Cgaps#1{}
\allowdisplaybreaks

\def\undersetbrace#1\to#2{\underbrace{#2}_{#1}}
\def\oversetbrace#1\to#2{\overbrace{#2}^{#1}}
\def\AMSunderset#1\to#2{\underset{#1}{#2}}
\def\AMSoverset#1\to#2{\overset{#1}{#2}}

\swapnumbers

\newtheorem*{prop*}{Proposition}

\newtheorem*{thm*}{Theorem}

\newtheorem*{lem*}{Lemma}

\newtheorem*{cor*}{Corollary}
\newenvironment{demo}[1]{\par\smallskip\noindent{\bf #1.}}{\par\smallskip}

\def\cit#1#2{\ifx#1!\cite{#2}\else#2\fi} 
\def\ign#1{}             
\def\o{\circ}
\def\X{\mathfrak X}
\def\al{\alpha}
\def\be{\beta}

\def\de{\delta}
\def\ep{\varepsilon}
\def\ze{\zeta}
\def\et{\eta}
\def\th{\theta}

\def\ka{\kappa}
\def\la{\lambda}
\def\rh{\rho}
\def\si{\sigma}

\def\ph{\varphi}
\def\ch{\chi}
\def\ps{\psi}
\def\om{\omega}
\def\Ga{\Gamma}
\def\De{\Delta}
\def\Th{\Theta}

\def\Ph{\Phi}

\def\i{^{-1}}
\def\x{\times}
\def\p{\partial}
\let\on=\operatorname

\def\AMSonly#1{}
\begin{document}

\title[Riemannian metrics on spaces
of curves, Hamiltonian approach]
{An overview of the Riemannian metrics on spaces
of curves using the Hamiltonian approach
}
\author{Peter W. Michor, David Mumford}
\address{
Peter W. Michor:
Fakult\"at f\"ur Mathematik, Universit\"at Wien,
Nordbergstrasse 15, A-1090 Wien, Austria; {\it and:}
Erwin Schr\"odinger Institut f\"ur Mathematische Physik,
Boltzmanngasse 9, A-1090 Wien, Austria
}
\email{Peter.Michor@univie.ac.at}
\address{
David Mumford:
Division of Applied Mathematics, Brown University,
Box F, Providence, RI 02912, USA}
\email{David\_{}Mumford@brown.edu}
\thanks{Both authors were supported by 
    NSF-Focused Research Group: 
    The geometry, mechanics, and statistics of the infinite dimensional
    shape manifolds.
    PWM was supported by FWF Project P~17108, and thanks the Program for
     Evolutionay Dynamics, Harvard University, for hospitality}
\date{\today}
\begin{abstract}
Here shape space is either the manifold of simple closed smooth
unparameterized curves in $\mathbb R^2$ or is the orbifold of immersions
from $S^1$ to $\mathbb R^2$ modulo the group of diffeomorphisms of $S^1$.
We investige several Riemannian metrics on shape space: $L^2$-metrics
weighted by expressions in length and curvature. These include a scale
invariant metric and a Wasserstein type metric which is sandwiched between
two length-weighted metrics. Sobolev metrics of order $n$ on curves are
described. Here the horizontal projection of a tangent field is given by a
pseudo-differential operator. Finally the metric induced from the Sobolev
metric on the group of diffeomorphisms on $\mathbb R^2$is treated. Although
the quotient metrics are all given by pseudo-differential operators, their
inverses are given by convolution with smooth kernels. We are able to prove
local existence and uniqueness of solution to the geodesic equation for
both kinds of Sobolev metrics.

We are interested in all conserved quantities, so the paper starts with the
Hamiltonian setting and computes conserved momenta and geodesics in general
on the space of immersions. For each metric we compute the geodesic
equation on shape space. In the end we sketch in some examples the
differences between these metrics.
\end{abstract}
\def\LaTeXonly{}

\maketitle

1. Introduction --- multiple Riemannian metrics on the space of curves
\newline\indent
2. The Hamiltonian approach
\newline\indent
3. Almost local Riemannian metrics on $\on{Imm}(S^1,\mathbb R^2)$ and on $B_i$.
\newline\indent
4. Sobolev metrics on $\on{Imm}(S^1,\mathbb R^2)$ and on $B_i$
\newline\indent
5. Sobolev metrics on $\on{Diff}(\mathbb R^2)$ and its quotients
\newline\indent
6. Examples

\section{Introduction --- multiple Riemannian metrics on the space of curves}
\label{nmb:1}

Both from a mathematical and a computer vision point of view, it is of
great interest to understand the space of simple closed curves in the
plane. Mathematically, this is arguably the simplest infinite-dimensional
truly nonlinear space. From a vision perspective, one needs to make human
judgements like `such-and-such shapes are similar, but such-and-such are
not' into precise statements. The common theory which links these two
points of view is the study of the various ways in which the space of
simple closed curves can be endowed with a Riemannian metric. From a vision
perspective, this converts the idea of similarity of two shapes into a
quantitative metric. From a mathematical perspective, a Riemannian metric
leads to geodesics, curvature and diffusion and, hopefully, to an
understanding of the global geometry of the space. Much work has been done
in this direction recently (see for example \cite{MennYezzi, MM1, MillerY,
MillerTY, Y}). The purpose of the present paper is two-fold. On the one
hand, we want to survey the spectrum of Riemannian metrics which have been
proposed (omitting, however, the Weil-Peterson metric). On the other hand,
we want to develop systematically the Hamiltonian approach to analyzing
these metrics.

Next, we define the spaces which we will study and introduce the notation we
will follow throughout this paper. To be precise, by a curve we mean a
$C^\infty$ simple closed curve in the plane. The space of these will be
denoted $B_e$. We will consider two approaches to working with this space.
In the first, we use parametrized curves and represent $B_e$ as the
quotient:
$$ B_e \cong \on{Emb}(S^1,\mathbb R^2)/\on{Diff}(S^1)$$
of the smooth Fr\'echet manifold of $C^\infty$ embeddings of $S^1$ in the
plane modulo the group of $C^\infty$ diffeomorphisms of $S^1$. In this
approach, it is natural to consider all possible {\it immersions} as well
as embeddings, and thus introduce the larger
space $B_i$ as the quotient of the space of $C^\infty$ immersions by the
group of diffeomorphisms of $S^1$:
\begin{equation*}
\begin{array}{ccccc}
\on{Emb}(S^1,\mathbb R^2) & \longrightarrow & \on{Emb}(S^1,\mathbb 
R^2)/\on{Diff}(S^1) & \cong & B_e\\
\cap & & \cap & & \cap\\
\on{Imm}(S^1,\mathbb R^2) & \longrightarrow & \on{Imm}(S^1,\mathbb 
R^2)/\on{Diff}(S^1) & \cong & B_i
\end{array}
\end{equation*}

In the second approach, we use the group of diffeomorphisms 
$\on{Diff}(\mathbb R^2)$ of the plane, where, more precisely, this is
either the group of all diffeomorphisms equal to the identity outside a
compact set or the group of all diffeomorphisms which decrease rapidly to
the identity. Let $\De$ be the unit circle in the plane. This group has two
subgroups, the normalizer and the centralizer of $\De$ in
$\on{Diff}(\mathbb R^2)$:
\begin{equation*}
\begin{array}{ccccc}
\on{Diff}^0(\mathbb R^2, \De) & \subset &\on{Diff}(\mathbb R^2, \De)&
\subset & \on{Diff}(\mathbb R^2) \\
\parallel & & \parallel & &\\
\left\{\ph \mid \ph|_\De \equiv \text{id}_\De \right\} & & \left\{\ph \mid
\ph(\De)= \De \right\}& &
\end{array}
\end{equation*}
Let $i \in \on{Emb}(S^1,\mathbb R^2)$ be the basepoint $i(\th) =
(\sin(\th), \cos(\th))$ carrying $S^1$ to the unit circle $\De$.
The group $\on{Diff}(\mathbb R^2)$ acts on the space $\on{Emb}(S^1,\mathbb
R^2)$ of embeddings by composition on the left. The action on the space of
embeddings is transitive (e.g., choose an isotopy between two embedded
circles, transform and extend its velocity field into a time-dependent
vector field with compact support on $\mathbb R^2$ and integrate it to a
diffeomorphism).
$\on{Diff}^0(\mathbb R^2, \De)$ is the subgroup which fixes the base
point $i$. Thus we can represent $\on{Emb}(S^1,\mathbb R^2)$ as the coset
space $\on{Diff}(\mathbb R^2)/\on{Diff}^0(\mathbb R^2,\De)$.

Furthermore $\on{Diff}^0(\mathbb R^2, \De)$ is a {\it normal} subgroup of
$\on{Diff}(\mathbb R^2, \De)$, and the quotient of one by the other is
nothing other than $\on{Diff}(\De)$, the diffeomorphism group of the unit
circle. So $\on{Diff}(\De)$ acts on the coset space $\on{Diff}(\mathbb
R^2)/\on{Diff}^0(\mathbb R^2,\De)$ with quotient the coset space
$\on{Diff}(\mathbb R^2)/\on{Diff}(\mathbb R^2,\De)$. Finally, under the
identification of $\on{Diff}(\mathbb R^2)/\on{Diff}^0(\mathbb R^2,\De)$
with $\on{Emb}(S^1,\mathbb R^2)$, this action is the same as the previously
defined one of $\on{Diff}(S^1)$ on $\on{Emb}(S^1,\mathbb R^2)$. This is
because if $c = \ph \circ i \in \on{Emb}(S^1,\mathbb R^2)$, and $\ps \in
\on{Diff}(\mathbb R^2, \De)$ satisfies $\ps(i(\th)) = i(h(\th)), h \in
\on{Diff}(S^1)$, then the action of $\ps$ carries $\ph$ to $\ph \circ \ps$
and hence $c$ to $\ph \circ \ps \circ i = \ph \circ i \circ h = c \circ h$.

{\it All } the spaces and maps we have introduced can be combined in one
commutative diagram:
\begin{equation*}
\begin{array}{ccccc}
\on{Diff}(\mathbb R^2) & & & & \\
\downarrow & & & &\\
\on{Diff}(\mathbb R^2) / \on{Diff}^0(\mathbb R^2, \De) &
\stackrel{\approx}{\longrightarrow} & \on{Emb}(S^1, \mathbb R^2) & \subset
& \on{Imm}(S^1, \mathbb R^2) \\
\downarrow & & \downarrow & & \downarrow \\
\on{Diff}(\mathbb R^2) / \on{Diff}(\mathbb R^2, \De) &
\stackrel{\approx}{\longrightarrow} & B_e & \subset & B_i
\end{array}\end{equation*}
See \cite{MM1} and \cite{KoM} for the homotopy type of the
spaces $\on{Imm}(S^1,\mathbb R^2)$ and $B_i$. 

What is the infinitesimal version of this? We will use the notation
$\X(\mathbb R^2)$ to denote the Lie algebra of $\on{Diff}(\mathbb R^2)$,
i.e., either the space of vector fields on $\mathbb R^2$ with compact
support or the space of rapidly decreasing vector fields. As for any Lie
group, the tangent bundle $T\on{Diff}(\mathbb R^2)$ is the product
$\on{Diff}(\mathbb R^2) \times \X(\mathbb R^2)$ by either right or left
multiplication. We choose right so that a tangent vector to
$\on{Diff}(\mathbb R^2)$ at $\ph$ is given by a vector field $X$
representing the infinitesimal curve $\ph \mapsto \ph(x,y) + \ep
X(\ph(x,y)$).

Fix $\ph \in \on{Diff}(\mathbb R^2)$ and let it map to $c = \ph \circ i \in
\on{Emb}(S^1,\mathbb R^2)$ and to the curve $C=\on{Im}(c) \subset \mathbb
R^2$ on the three levels of the above diagram. A tangent vector to
$\on{Emb}(S^1, \mathbb R^2)$ at $c$ is given by a vector field $Y$ to
$\mathbb R^2$ along the map $c$, and the vertical map of tangent vectors
simply takes the vector field $X$ defined on all of $\mathbb R^2$ and
restricts it to the map $c$, i.e.\ it takes the values $Y(\th) =
X(c(\th))$. Note that if $c$ is an embedding, a vector field along $c$ is
the same as a vector field on its image $C$. A tangent vector to $B_e$ at
the image curve $C$ is given by a vector field $Y$ along $C$ {\it modulo}
vector fields tangent to $C$ itself. The vertical map on tangent vectors
just takes the vector field $X$ along $c$ and treats it modulo vector
fields tangent to $c$. However, it is convenient to represent a tangent
vector to $B_e$ or $B_i$ at $C$ not as an equivalence class of vector
fields along $C$ but by their unique representative which is everywhere
{\it normal} to the curve $C$. This makes $T_CB_i$ the space of all
normal vector fields to $C \subset \mathbb R^2$.

In both approaches, we will put a Riemannian metric on the top 
space, i.e. $\on{Imm}(S^1,\mathbb R^2)$ or $\on{Diff}(\mathbb R^2)$, 
which makes the map to the quotient $B_i$ or to a coset space of
$\on{Diff}(\mathbb R^2)$ into a {\it Riemannian submersion}. In general,
given a diffeomorphism $f:A \rightarrow B$ with a surjective tangent map
and a metric $G_a(h,k)$ on $A$, $f$ is a submersion if it has the following
property: first split the tangent bundle to $A$ into the subbundle
$TA^\top$ tangent to the fibres of $f$ and its perpendicular
$TA^\bot$ with respect to $G$ (called the horizontal bundle).
Then, under the isomorphisms
$df:TA^\bot_a \stackrel{\approx}{\rightarrow} TB_{f(a)}$, the restriction
of the $A$-metric to the horizontal subbundle is required to define a
metric on $TB_b$, independent of the choice of the point $a \in f\i(b)$ in
the fiber. In this
way we will define Riemannian metrics on all the spaces in our diagram
above. Submersions have a very nice effect on geodesics: the geodesics on
the quotient space $B$ are exactly the images of the geodesics on the top
space $A$ which are perpendicular at one, and hence at all, points to the
fibres of the map $f$ (or, equivalently, their tangents are in the
horizontal subbundle).

On $\on{Diff}(\mathbb R^2)$, we will consider only right invariant metrics.
These are given by putting a positive definite inner product $G(X,Y)$ on
the vector space of vector fields to $\mathbb R^2$, and translating this to
the tangent space above each diffeomorphism $\ph$ as above. That is, the
length of the infinitesimal curve $\ph \mapsto \ph + \ep X\circ \ph$ is
$\sqrt{G(X,X)}$. Then the map from $\on{Diff}(\mathbb R^2)$ to any of its
right coset spaces will be a Riemannian submersion, hence we get metrics on
all these coset spaces.

A Riemannian metric on $\on{Imm}(S^1,\mathbb R^2)$ is just a family of
positive definite inner products $G_c(h,k)$ where $c \in
\on{Imm}(S^1,\mathbb R^2)$ and $h,k \in C^\infty(S^1,\mathbb R^2)$
represent vector fields on $\mathbb R^2$ along $c$. We require that our
metrics will be invariant under the action of $\on{Diff}(S^1)$, hence the
map dividing by this action will be a Riemannian submersion. Thus we will
get Riemannian metrics on $B_i$: these are given by a family of inner
products as above such that $G_c(h,k) \equiv 0$ if $h$ is tangent to $c$,
i.e., $\langle h(\th), c_\th(\th) \rangle \equiv 0$ where $c_\th:=\p_\th c$.

When dealing with parametrized curves or, more generally, immersions, we
will use the following terminology. Firstly, the immersion itself is
usually denoted by:
$$ c(\th):S^1 \rightarrow \mathbb R^2$$
or, when there is a family of such immersions:
$$ c(\th,t):S^1\times I \rightarrow \mathbb R^2.$$
The parametrization being usually irrelevant, we work mostly with 
arclength $ds$, arclength derivative $D_s$ and the unit tangent 
vector $v$ to the curve:
\begin{align*} ds &= |c_\th| d\th \\
D_s &= \partial_\th/|c_\th| \\
v &= c_\th/|c_\th|
\end{align*}
An important caution is that when you have a family of curves 
$c(\th,t)$, then $\partial_\th$ and $\partial_t$ commute but $D_s$ 
and $\partial_t$ don't because $|c_\th|$ may have a $t$-derivative. 
Rotation through 90 degrees will be denoted by:
$$ J = \left( \begin{array}{cc} 0 & -1 \\ 1 & 0 \end{array} \right).$$
The unit normal vector to the image curve is thus
$$n = Jv.$$
Thus a Riemannian metric on $B_e$ or $B_i$ is given by inner products
$G_C(a,b)$ where $a.n$ and $b.n$ are any two vector fields along $C$ normal
to $C$ and $a,b \in C^\infty(C,\mathbb R)$. Another important piece of
notation that we will use concerns directional derivatives of functions
which depend on several variables. Given a function $f(x,y)$ for instance,
we will write:
$$ D_{(x,h)}f \text{ or } df(x)(h) \text{ as shorthand for } \partial_t|_0 f(x+th,y).$$
Here the $x$ in the subscript will indicate which variable is changing and
the second argument $h$ indicates the direction. This applies even if one
of the variables is a curve $C \in B_i$ and $h$ is a normal vector field.

The simplest inner product on the tangent bundle to $\on{Imm}(S^1,\mathbb
R^2)$ is:
$$ G^0_c(h,k) = \int_{S^1} \langle h(\th), k(\th) \rangle \cdot ds.$$
Since the differential $ds$ is invariant under the action of the 
group $\on{Diff}(S^1)$, the map to the quotient $B_i$ is a Riemannian 
submersion for this metric. A tangent vector $h$ to 
$\on{Imm}(S^1,\mathbb R^2)$ is perpendicular to the orbits of 
$\on{Diff}(S^1)$ if and only if $\langle h(\th), v(\th) \rangle 
\equiv 0$, i.e. $h$ is a multiple $a.n$ of the unit normal. This is the
same subbundle as above, so that, for this metric, the {\it horizontal}
subspace of the tangent space is the natural splitting. Finally, the
quotient metric is given by
$$ G^0_c(a\cdot n, b\cdot n) = \int_{S^1} a.b.ds.$$

All the metrics we will look at will be of the form:
$$ G_c(h,k) = \int_{S^1} \langle Lh, k \rangle \cdot ds$$
where $L$ is a positive definite operator on the vector-valued functions
$h:S^1 \rightarrow \mathbb R^2$. The simplest such $L$ is simply
multiplication by some function $\Ph_c(\th)$. However, it will turn out
that most of the metrics involve $L$'s which are differential or {\it
pseudo-differential} operators. For these, the horizontal subspace is not
the natural splitting, so the quotient metric on $B_e$ and $B_i$ involves
restricting $G_c$ to different sub-bundles and this makes these operators
somewhat complicated. In fact, it is not guaranteed that the horizontal
subspace is spanned by $C^\infty$ vectors (in the sense that the full
$C^\infty$ tangent space is the direct sum of the vertical subspace and the
horizontal $C^\infty$ vectors). When dealing with metrics on
$\on{Diff}(\mathbb R^2)$ and vertical subspaces defined by the subgroups
above, this does happen. In this case, the horizontal subspace must be
taken using less smooth vectors. 

In all our cases, $L\i$ will be a simpler operator than $L$: this is
because the tangent spaces to $B_e$ or $B_i$ are {\it quotients} of the
tangent spaces to the top spaces $\on{Diff}$ or $\on{Imm}$ where the
metrics are most simply defined, whereas the cotangent spaces to $B_e$ or
$B_i$ are subspaces of the cotangent spaces of the space `above'. The dual
inner product on the cotangent space is given by the inverse operator $L\i$
and in all our cases this will be an integral operator with a simple
explicit kernel. A final point: we will use a constant $A$ when terms with
different physical `dimensions' are being added in the operator $L$. Then
$A$ plays the role of fixing a scale relative to which different geometric
phenomena can be expected.

Let us now describe in some detail the contents of this paper and the
metrics. First, in section \ref{hamiltoniansection}, we introduce the
general Hamiltonian formalism. This is, unfortunately, more technical than
the rest of the paper. First we consider general Riemannian metrics on the
space of immersions which admit Christoffel symbols. We express this as the
existence of two kinds of gradients. Since the energy function is not even
defined on the whole cotangent bundle of the tangent bundle we pull back to
the tangent bundle the canonical symplectic structure on the cotangent
bundle. Then we determine the Hamiltonian vector field mapping and, as a
special case, the geodesic equation. We detemine the equivariant moment
mapping for several group actions on the space of immersions: the action of
the reparametrization group
$\on{Diff}(S^1)$, of the motion group of $\mathbb R^2$, and also of the
scaling group (if the metric is scale invariant). Finally the invariant
momentum mapping on the group $\on{Diff}(\mathbb R^2)$ is described.

Section \ref{almostlocalmetrics} is then devoted to applying the Hamiltonian
procedure to {\it almost local metrics}: these are the metrics in which $L$
is multiplication by some function $\Ph$. Let $\ell_c=\int_{S^1} ds $ be
the length of the curve $c$ and let
$$ \ka_c(\th)= \langle  n(\th),D_s(v)(\th)\rangle = \langle Jc_\th, 
c_{\th \th} \rangle/ |c_\th|^3$$
be the curvature of $c$ at $c(\th)$. Then for any auxiliary function 
$\Phi(\ell,\ka)$, we can define a weighted Riemannian metric:
$$ G^\Phi_c(h,k) = \int_{S^1} \Phi(\ell_c, \ka_c(\th)) \cdot h(\th) 
k(\th)\cdot ds.$$
The motivation for introducing weights is simply that, for any 2 curves in
$B_i$, the infimum of path lengths in the $G^0$ metric for paths joining
them is zero,
see \cite{MM1, MM2}. For all these metrics, the horizontal subspace is
again the set of tangent vectors $a(\th) n(\th)$, so the metric on $B_i$ is
simply
$$ G^\Phi_c(a\cdot n, b\cdot n) = \int_{S^1} \Phi(\ell_c, \ka_c(\th))
a(\th)b(\th)\cdot ds.$$
We will determine the geodesic equation, the momenta and the sectional
curvature for all these metrics. The formula for sectional curvature is
rather complicated but for special $\Ph$, it is quite usable.

We will look at several special cases. The weights
$$ \Phi(\ell,\ka) = 1 + A \ka^2$$
were introduced and studied in \cite{MM1}. As we shall see, this metric is
also closely connected to the Wasserstein metric on probability measures
(see \cite{ambrosio}), if we assign to a curve $C$ the probability measure
given by scaled arc length. We show that it is sandwiched between the conformal
metric $G^{\ell^{-1}}$ and $G^{\Ph_W}$ where
$\Ph_W=\ell^{-1}+\frac1{12}\ell\ka^2$. Weights of the form
$$ \Phi(\ell,\ka) = f(\ell)$$
were studied in \cite{MennYezzi} and independently by \cite{Shah}. The latter are
attractive because they give metrics which are conformally equivalent to
$G^0$. These metrics are a borderline case between really stable metrics on
$B_e$ and the metric $G^0$ for which path length shrinks to 0: for them,
the infimum of path lengths is positive but at least some paths seem to
oscillate wildly when their length approaches this infimum. Another very
interesting case is:
$$ \Phi(\ell, \ka) =\ell^{-3} + A|\ka|^2\ell^{-1}$$
because this metric is scale-invariant. 

A more standard approach to strengthening $G^0$ is to introduce higher
derivatives. In section \ref{G^{imm,n}}, we follow the Sobolev approach
which puts a metric on $\on{Imm}(S^1, \mathbb R^2)$ by:
\begin{align*}
G^{\text{imm},n}_c(h,k) &= \int_{S^1} \sum^n_{i=0} \langle D^i_s h, D^i_s k
\rangle ds \\
&= \int_{S^1} \langle Lh,k \rangle ds, \quad \text{where }L = \sum_{i=0}^n
(-1)^i D_s^{2i}
\end{align*}
However, the formulas we get are substantially simpler and $L\i$ has an
elegant expression is we take the equivalent metric:
\begin{align*}
G^{\text{imm},n}_c(h,k) &= \int_{S^1} \left( \langle h,k \rangle + A.\langle D^n_s h,
D^n_s k
\rangle \right) ds \\
&= \int_{S^1} \langle Lh,k \rangle ds, \quad \text{where }L = I + (-1)^n A
D_s^{2n}
\end{align*}
We apply the Hamiltonian procedure to this metric. Here the horizontal
space of all vectors in the tangent space $T\on{Imm}(S^1,\mathbb R^2)$
which are
$G^{\text{imm},n}$-orthogonal to the reparametrization orbits, is very
different from the natural splitting in \S \ref{almostlocalmetrics}. The
decomposition of a vector into horizontal and vertical parts involves pseudo
differential operators, and thus also the horizontal geodesic equation is
an integro-differential equation. However, its inverse $L\i$ is an integral
operator whose kernel has a simple expression in terms of arc length
distance between 2 points on the curve and their unit normal vectors.

For this metric, we work out the geodesic equation and prove that the
geodesic flow is well posed in the sense that we have local existence and
uniqueness of solutions in $\on{Imm}(S^1,\mathbb R^2)$ and in $B_i$.
Finally we discuss a little bit a scale invariant version of the metric
$G^{\text{imm},n}$. For the simplest of these metrics, the scaling
invariant momentum along a geodesic turns out to be the time derivative of
$\log(\ell)$. At this time, we do not know the sectional curvature for this
metric.

In the next section \ref{G^{diff,n}}, we start with the basic right
invariant metrics on $\on{Diff}(\mathbb R^2)$ which are given by the
Sobolev $H^n$-inner product on $\X(\mathbb R^2)$.
\begin{align*}
H^n(X,Y) &=\sum_{i,j \ge 0, i+j\le n} \frac{A^{i+j}n!}{i!j!(n-i-j)!}
\iint_{\mathbb R^2}
   \langle \p_{x}^i\p_{y}^j X,\p_{x}^i\p_{y}^j Y \rangle\,dx^1 dx^2 \\
&= \iint_{\mathbb R^2} \langle LX,Y \rangle dx.dy, \quad \text{where }
L=(1-A \De)^n, \De = \p_x^2 + \p_y^2.
\end{align*}
These metrics have been extensively studied by Miller, Younes and Trouv\'e
and their collaborators \cite{BegMillerTY, MillerY, MillerTY, YT}. Since
these metrics are right invariant, all maps to coset spaces
$\on{Diff}(\mathbb R^2) \rightarrow \on{Diff}(\mathbb R^2)/H$ are
submersions. In particular, this metric gives a quotient metric on
$\on{Emb}(S^1, \mathbb R^2)$ and $B_e$ which we will denote by
$G^{\text{diff},n}_c(h,k)$. In this case, the inverse $L\i$ of the operator
defining the metric is an integral operator with a kernel given by a
classical Bessel function applied to the distance in $\mathbb R^2$
between 2 points on the curve. We will derive the geodesic equations: they
are all in the same family as fluid flow equations. We prove well posedness
of the geodesic equation on $\on{Emb}(S^1,\mathbb R^2)$ and on $B_e$.
Although there is a formula of Arnold \cite{Arnold} for the sectional
curvature of any right-invariant metric on a Lie group, we have not
computed sectional curvatures for the quotient spaces.

In the final section \ref{examples}, we study two examples to make clear
the differences between the various metrics. The first example is the
geodesic formed by the set of all concentric circles with fixed center. We
will see how this geodesic is complete when the metric is reasonably
strong, but incomplete in most `borderline' cases. The second example takes
a fixed `cigar-shaped' curve $C$ and compares the unit balls in the tangent
space $T_C B_e$ given by the different metrics.

\section{The Hamiltonian approach}\label{hamiltoniansection}

In our previous papers, we have derived the geodesic equation in our various 
metrics by setting the first variation of the energy of a path equal to 0. 
Alternately, the geodesic equation is the Hamiltonian flow associated 
to the first fundamental form (i.e. the length-squared function given 
by the metric on the tangent bundle). 
The Hamiltonian approach also provides a mechanism for converting symmetries 
of the underlying Riemannian manifold into conserved quantities, the momenta. 
We first need to be quite formal and lay out the basic definitions, 
esp.\ distinguishing between the tangent and cotangent bundles rather
carefully: The former consists of smooth vector fields along immersions
whereas the latter is comprised of 1-currents along immersions. Because of
this we work on the tangent bundle and we pull back the symplectic form
from the contantent bundle to $T\on{Imm}(S^1,\mathbb R^2)$. We use the
basics of symplectic geometry and momentum mappings on cotangent bundles in
infinite dimensions, and we explain each step. 
See \cite{MGeomEvol},~section~2, for a detailed exposition in similar
notation as used here.

\subsection{The setting}
\label{setting}
Consider as above the smooth Fr\'echet manifold $\on{Imm}(S^1,\mathbb R^2)$
of all immersions $S^1\to \mathbb R^2$ which is an open subset of
$C^\infty(S^1,\mathbb R^2)$. The tangent bundle is $T\on{Imm}(S^1,\mathbb
R^2)=\on{Imm}(S^1,\mathbb R^2)\x C^\infty(S^1,\mathbb R^2)$,
and the cotangent bundle is $T^*\on{Imm}(S^1,\mathbb
R^2)=\on{Imm}(S^1,\mathbb R^2)\x \mathcal D(S^1)^2$ where the second factor
consists of pairs of periodic distributions.

We consider smooth Riemannian metrics on $\on{Imm}(S^1,\mathbb R^2)$, i.e.,
smooth mappings
\begin{align*}
&G:\on{Imm}(S^1,\mathbb R^2)\x C^\infty(S^1,\mathbb R^2)
  \x C^\infty(S^1,\mathbb R^2)\to \mathbb R \\&
(c,h,k)\mapsto G_c(h,k),\quad \text{ bilinear in }h,k
\\& G_c(h,h)>0 \quad\text{  for }h\ne0.
\end{align*}
Each such metric is {\it weak} in the sense that $G_c$, viewed as bounded
linear mapping
\begin{align*}
G_c:T_c\on{Imm}(S^1,\mathbb R^2)=C^\infty(S^1,\mathbb R^2)&\to
T_c^*\on{Imm}(S^1,\mathbb R^2)=\mathcal D(S^1)^2
\\
G: T\on{Imm}(S^1,\mathbb R^2)&\to T^*\on{Imm}(S^1,\mathbb R^2)
\\
G(c,h)&=(c,G_c(h,\;.\;))
\end{align*}
is injective, but can never be surjective. We shall need also its tangent
mapping 
\begin{align*}
TG:T(T\on{Imm}(S^1,\mathbb R^2))&\to T(T^*\on{Imm}(S^1,\mathbb R^2))
\end{align*}
We write a tangent vector to $T\on{Imm}(S^1,\mathbb R^2)$ in the form 
$(c,h;k,\ell)$ where $(c,h)\in T\on{Imm}(S^1,\mathbb R^2)$ is its foot point, 
$k$ is its vector component in the $\on{Imm}(S^1,\mathbb R^2)$-direction
and where $\ell$ is
its component in the $C^\infty(S^1,\mathbb R^2)$-direction. 
Then $TG$ is given by
$$
TG(c,h;k,\ell) = (c, G_c(h,\;.\;);k, D_{(c,k)}G_c(h,\;.\;)+G_c(\ell,\;.\;))
$$
Moreover, if
$X=(c,h;k,\ell)$ then we will write $X_1=k$ for its first vector component
and $X_2=\ell$ for the second vector component. 
Note that only these smooth functions on $\on{Imm}(S^1,\mathbb R^2)$ whose
derivative lies in the image of $G$ in the cotangent bundle have
$G$-gradients. This requirement has only to be satisfied for the first
derivative, for the higher ones it follows (see \cite{KM}). We shall denote
by $C^\infty_G(\on{Imm}(S^1,\mathbb R^2))$ the space of such smooth
functions.

We shall always assume that $G$ is invariant under the reparametrization group
$\on{Diff}(S^1)$, hence each such metric induces a Riemann-metric on the
quotient space $B_i(S^1,\mathbb R^2)=\on{Imm}(S^1,\mathbb
R^2)/\on{Diff}(S^1)$.

In the sequel we shall further assume that that {\it the weak Riemannian
metric $G$ itself admits $G$-gradients with respect to the variable $c$ in
the following sense:}
\begin{align*}
&\boxed{
D_{c,m}G_c(h,k) = G_c(m,H_c(h,k)) = G_c(K_c(m,h),k)
}\quad\text{  where }
\\
&H,K:\on{Imm}(S^1,\mathbb R^2)\x C^\infty(S^1,\mathbb R^2)
  \x C^\infty(S^1,\mathbb R^2)\to C^\infty(S^1,\mathbb R^2)
\\& \hspace*{2.5in}(c,h,k)\mapsto H_c(h,k),K_c(h,k)
\\&
\text{are smooth and bilinear in }h,k.
\end{align*}
Note that $H$ and $K$ could be expressed in (abstract) index notation as 
$g_{ij,k}g^{kl}$ and $g_{ij,k}g^{il}$.
We will check and compute these gradients for several concrete metrics below.

\subsection{The fundamental symplectic form on $T\on{Imm}(S^1,\mathbb R^2)$
induced by a weak Riemannian metric}
\label{symplectic form}
The basis of Hamiltonian theory is the natural 1-form on the cotangent
bundle $T^*\on{Imm}(S^1,\mathbb R^2)$ 
given by:
\begin{gather*}
\Th: T(T^*\on{Imm}(S^1,\mathbb R^2))=\on{Imm}(S^1,\mathbb R^2)\x
\mathcal D(S^1)^2\x C^\infty(S^1,\mathbb R^2)\x \mathcal D(S^1)^2 \to \mathbb R
\\
(c,\al;h,\be) \mapsto \langle  \al,h\rangle.
\end{gather*}
The pullback via the mapping $G: T\on{Imm}(S^1,\mathbb R^2)\to
T^*\on{Imm}(S^1,\mathbb R^2)$ of the 1-form $\Th$ is then:
$$ (G^*\Th)_{(c,h)}(c,h;k,\ell)=G_c(h,k).$$
Thus the symplectic form $\om=-dG^*\Th$ on $ T\on{Imm}(S^1,\mathbb R^2)$
can be computed as follows, where we use the constant vector fields
$(c,h)\mapsto (c,h;k,\ell)$:
\begin{align*}
\om_{(c,h)}&((k_1,\ell_1),(k_2,\ell_2)) =
-d(G^*\Th)((k_1,\ell_1),(k_2,\ell_2))|_{(c,h)} \\&
=-D_{(c,k_1)}G_c(h,k_2)- G_c(\ell_1,k_2)
+D_{(c,k_2)}G_c(h,k_1)+ G_c(\ell_2,k_1) 
\\&
=G_c\big(k_2,H_c(h,k_1)-K_c(k_1,h)\big) + G_c(\ell_2,k_1) - G_c(\ell_1,k_2)
\tag1
\end{align*}

\subsection{The Hamiltonian vector field mapping}
\label{hamiltonian}
Here we compute the Hamiltonian vectorfield $\on{grad}^\om(f)$ associated
to a smooth function $f$ on the tangent space $T\on{Imm}(S^1,\mathbb R^2)$, that is $f \in C^\infty_G(\on{Imm}(S^1,\mathbb R^2)\x C^\infty(S^1,\mathbb R^2))$
assuming that it has smooth $G$-gradients in both factors. See
\cite{KM},~section~48. Using the explicit formulas in \ref{symplectic
form}, we have:
\begin{align*}
&\om_{(c,h)}\left(\on{grad}^\om(f)(c,h),(k,\ell) \right) 
=\om_{(c,h)}\left((\on{grad}_1^\om(f)(c,h),\on{grad}_2^\om(f)(c,h)),(k,\ell) \right) =
\\&
= G_c\big(k, H_c\big(h,\on{grad}_1^\om(f)(c,h)\big)\big)
-G_c(K_c(\on{grad}_1^\om(f)(c,h),h),k)
\\&\quad
 + G_c(\ell,\on{grad}_1^\om(f)(c,h)) 
 - G_c(\on{grad}_2^\om(f)(c,h),k) 
\end{align*}
On the other hand, by the definition of the $\om$-gradient we have
\begin{align*}
&\om_{(c,h)}\left(\on{grad}^\om(f)(c,h),(k,\ell) \right) 
=df(c,h)(k,\ell)= D_{(c,k)}f(c,h) + D_{(h,\ell)}f(c,h)
\\&
=G_c(\on{grad}_1^{G}(f)(c,h),k) 
+G_c(\on{grad}_2^{G}(f)(c,h),\ell) 
\end{align*}
and we get the expression of the Hamiltonian vectorfield:
$$
\boxed{
\begin{aligned}
\on{grad}_1^\om(f)(c,h) &= \on{grad}_2^{G}(f)(c,h)\\
\on{grad}_2^\om(f)(c,h) &= \!-\! \on{grad}_1^{G}(f)(c,h)
\!+\! H_c\big(h,\on{grad}_2^{G}(f)(c,h)\big)
\!-\! K_c(\on{grad}_2^G(f)(c,h),h)
\end{aligned}
}
$$
Note that for a smooth function $f$ on $T\on{Imm}(S^1,\mathbb R^2)$ the
$\om$-gradient exists if and only if both $G$-gradients exist.

\subsection{The geodesic equation}
\label{geodesic}
The geodesic flow is defined by a vector field on $T\on{Imm}(S^1,\mathbb
R^2)$. One way to define this vector field is as the Hamiltonian vector
field of the energy function
$$ E(c,h)=\frac12 G_c(h,h),\qquad
E:\on{Imm}(S^1,\mathbb R^2)\x C^\infty(S^1,\mathbb R^2)\to \mathbb R .$$
The two partial $G$-gradients are:
\begin{align*}
G_c(\on{grad}_2^{G}(E)(c,h),\ell) &= d_2E(c,h)(\ell) = G_c(h,\ell)
\\
\on{grad}_2^{G}(E)(c,h)&=h
\\
G_c(\on{grad}_1^{G}(E)(c,h),k) &= d_1E(c,h)(k) 
= \tfrac12 D_{(c,k)}G_c(h,h) 
\\&
=\tfrac12 G_c(k,H_c(h,h))
\\
\on{grad}_1^{G}(E)(c,h) &= \tfrac12 H_c(h,h). 
\end{align*}
Thus the geodesic vector field is 
\begin{align*}
\on{grad}_1^\om(E)(c,h) &= h
\\
\on{grad}_2^\om(E)(c,h) &=  \tfrac12 H_c(h,h) - K_c(h,h)
\end{align*}
and the geodesic equation becomes:
\begin{align*}
&\begin{cases}
c_t &= h
\\
h_t &=  \tfrac12 H_c(h,h) - K_c(h,h)
\end{cases}
\quad\text{  or }\quad
\boxed{
c_{tt}= \tfrac12 H_c(c_t,c_t) - K_c(c_t,c_t)
}
\end{align*}
This is nothing but the usual formula for the geodesic flow using the
Christoffel symbols expanded out using the first derivatives of the metric
tensor.

\subsection{The momentum mapping for a $G$-isometric group action}
\label{momentum}
We consider now a (possibly infinite dimensional regular) Lie group with
Lie algebra $\mathfrak g$ with a right action $g\mapsto r^g$ by isometries
on $\on{Imm}(S^1,\mathbb R^2)$. If $\X(\on{Imm}(S^1,\mathbb R^2))$ denotes
the set of vector fields on $\on{Imm}(S^1,\mathbb R^2)$, we can specify
this action by the fundamental vector field mapping $\ze:\mathfrak g\to
\X(\on{Imm}(S^1,\mathbb R^2))$, which will be a bounded Lie algebra
homomorphism. The fundamental vector field $\ze_X, X \in \mathfrak g$ is
the infinitesimal action in the sense:
$$\ze_X(c)=\p_t|_0 r^{\exp(tX)}(c).$$ 
We also consider the tangent prolongation of this action on
$T\on{Imm}(S^1,\mathbb R^2)$ where the fundamental vector field is given by
\begin{equation*}
\ze_X^{T\on{Imm}}:(c,h)\mapsto (c,h;\ze_X(c),D_{(c,h)}(\ze_X)(c)=:\ze'_X(c,h))
\end{equation*}
The basic assumption is that the action is by isometries, 
$$ G_c(h,k)=((r^g)^*G)_c(h,k)=G_{r^g(c)}(T_c(r^g)h,T_c(r^g)k). $$
Differentiating this equation at $g=e$ in the direction $X\in \mathfrak g$
we get
\begin{equation*}
0=D_{(c,\ze_X(c))}G_c(h,k)+G_c(\ze'_X(c,h),k)+G_c(h,\ze'_X(c,k))
\tag{1}\end{equation*}
The key to the Hamiltonian approach is to define the group action by
Hamiltonian flows. To do this, we define the {\it momentum map} 
$j:\mathfrak g\to C^\infty_G(T\on{Imm}(S^1,\mathbb R^2), \mathbb R)$ by:
$$\boxed{j_X(c,h) = G_c(\ze_X(c),h).}$$
Equivalently, since this map is linear, it is often written as a map
$$\mathcal J:T\on{Imm}(S^1,\mathbb R^2) \rightarrow \mathfrak g',\qquad
  \langle \mathcal J(c,h),X \rangle = j_X(c,h).$$
The main property of the momentum map is that it fits into the following
commuta\-ti\-ve diagram and is a homomorphism of Lie algebras:
\begin{displaymath}
\cgaps{0.8;0.7;0.7;0.8}\xymatrix{
H^0(T\on{Imm}) \ar[r]^{i} &  
     C^\infty_G(T\on{Imm},\mathbb R)  \ar[rr]^{\on{grad}^{\om}} & &  
     \X(T\on{Imm},\om) \ar[r] & H^1(T\on{Imm})
\\
 & & \mathfrak g \ar[lu]^{j} \ar[ru]_{\ze^{T\on{Imm}}} & &
}\end{displaymath}
where $\X(T\on{Imm},\om)$ is the space of vector fields on $T\on{Imm}$
whose flow leaves $\om$ fixed. We need to check that:
\begin{align*}
\ze_X(c)=\on{grad}^\om_1(j_X)(c,h) &= \on{grad}^{G}_2(j_X)(c,h) \\
\ze'_X(c,h)=\on{grad}_2^\om(j_X)(c,h) &= - \on{grad}_1^{G}(j_X)(c,h) +
H_c(h,\ze_X(c)) -K_c(\ze_X(c),h)
\end{align*}
The first equation is obvious. To verify the second equation, we take its
inner product with some $k$ and use:
\begin{align*}
G\big(k,\on{grad}_1^{G}(j_X)(c,h)\big) &= D_{(c,k)}j_X(c,h) =
D_{(c,k)}G_c(\ze_X(c),h) + G_c(\ze'_X(c,k),h)
\\&
= G_c(k,H_c(\ze_X(c),h)) + G_c(\ze'_X(c,k),h).
\end{align*}
Combining this with \thetag{1}, the second equation follows. Let us check
that it is also a homomorphism of Lie algebras using the Poisson bracket:
\begin{align*}
\{j_X,j_Y\}(c,h)&=d j_Y(c,h)(\on{grad}^\om_1(j_X)(c,h),\on{grad}^\om_2(j_X)(c,h))
\\&
=d j_Y(c,h)(\ze_X(c),\ze_X'(c,h))
\\&
=D_{(c,\ze_X(c))}G_c(\ze_Y(c),h) + G_c(\ze_Y'(c,\ze_X(x)),h) +
G_c(\ze_Y(c),\ze_X'(c,h))
\\&
= G_c(\ze_Y'(c,\ze_X(c))-\ze_X'(c,\ze_Y(c)),h) \qquad\text{  by \thetag{1}}
\\&
= G_c([\ze_X,\ze_Y](c),h) = G_c(\ze_{[X,Y]}(c),h) = j_{[X,Y]}(c).
\end{align*}

Note also that $\mathcal J$ is equivariant for the group action,
by the following arguments: For $g$ in the Lie group let $r^g$ be the right
action on $\on{Imm}(S^1,\mathbb R^2)$, then $T(r^g)\o \ze_X\o (r^g)\i =
\ze_{\on{Ad}(g\i)X}$. Since $r^g$ is an isometry the mapping $T(r^g)$ is a
symplecto\-mor\-phism for $\om$, thus $\on{grad}^\om$ is equivariant. Thus
$j_X\o T(r^g)= j_{\on{Ad}(g)X}$ plus a possible constant which we can rule
out since $j_X(c,h)$ is linear in $h$.

By Emmy Noether's theorem, along any geodesic $t\mapsto c(t,\;.\;)$ this
momentum mapping is constant, thus for any $X\in\mathfrak g$ we have
\begin{displaymath}
\boxed{\quad
\langle \mathcal J(c,c_t),X \rangle = j_X(c,c_t) = G_c(\ze_X(c),c_t) \quad\text{  is
constant in }t.
\quad}
\end{displaymath}

We can apply this construction to the following group actions on
$\on{Imm}(S^1,\mathbb R^2)$.
\begin{itemize}
\item
The smooth right action of the group
$\on{Diff}(S^1)$ on $\on{Imm}(S^1,\mathbb R^2)$, 
given by composition from the right:
$c\mapsto c\o \ph$ for $\ph\in\on{Diff}(S^1)$. 
For $X\in\X(S^1)$ the fundamental vector field is then given by
$$ \ze^{\on{Diff}}_X(c) = \ze_X(c) = \p_t|_0 (c\o \on{Fl}^X_t) = c_\th.X $$
where $\on{Fl}^X_t$ denotes the flow of $X$.
The {\it reparametrization momentum}, for any vector field $X$ on $S^1$ is thus:
$$ j_X(c,h) = G_c(c_\th.X, h).$$
Assuming the metric is reparametrization invariant, it follows that on any
geodesic $c(\th,t)$, the expression $G_c(c_\th.X,c_t)$ is constant for all $X$.
\item
The left action of the Euclidean motion group 
$M(2)=\mathbb R^2\rtimes SO(2)$ on $\on{Imm}(S^1,\mathbb R^2)$ given by
$c\mapsto e^{aJ}c+ B$ for $(B,e^{aJ})\in \mathbb R^2\x SO(2)$. The 
fundamental vector field mapping is
\begin{align*}
\ze_{(B,a)}(c)&= aJc+B
\end{align*}
The {\it linear momentum} is thus $G_c(B,h), B \in \mathbb R^2$ and if the
metric is trans\-la\-tion invariant, $G_c(B,c_t)$ will be constant along
geodesics. The {\it angular momentum} is similarly $G_c(Jc,h)$ and if the
metric is rotation invariant, then $G_c(Jc,c_t)$ will be constant along
geodesics.

\item 
The action of the scaling group of $\mathbb R$ given by $c\mapsto e^r c$,
with fundamental vector field $\ze_a(c)=a.c$. 
If the metric is scale invariant, then
the {\it scaling momentum} $G_c(c,c_t)$ will also be invariant along
geodesics.
\end{itemize}

\subsection{Metrics and momenta on the group of diffeomorphisms}
Very similar things happen when we consider metrics on the group
$\on{Diff}(\mathbb R^2)$. As above, the tangent space to $\on{Diff}(\mathbb
R^2)$ at the identity is the vector space of vector fields $\X(\mathbb
R^2)$ on $\mathbb R^2$ and we can identify $T\on{Diff}(\mathbb R^2)$ with
the product $\on{Diff}(\mathbb R^2) \times \X(\mathbb R^2)$ using right
multiplication in the group to identify the tangent at a point $\ph$ with
that at the identity. The definition of this product decomposition means
that right multiplication by $\ps$ carries $(\ph,X)$ to $(\ph \circ \ps,
X)$. As usual, suppose that conjugation $\ph \mapsto \ps \circ \ph \circ
\ps\i$ has the derivative at the identity given by the linear operator
$\on{Ad}_\ps$ on the Lie algebra $\X(\mathbb R^2)$. It is easy to calculate
the explicit formula for Ad:
$$ \on{Ad}_\ps(X) = (D\ps \cdot X)\circ \ps\i.$$
Then left multiplication by $\ps$ on $\on{Diff}(\mathbb R^2) \times
\X(\mathbb R^2)$ is given by $(\ph,X) \mapsto (\ps \circ \ph,
\on{Ad}_\ps(X))$.
We now want to carry over the ideas of \ref{momentum} replacing the space
$\on{Imm}(S^1,\mathbb R^2)$ by $\on{Diff}(\mathbb R^2)$ and the group
action there by the right action of $\on{Diff}(\mathbb R^2)$ on itself. The
Lie algebra $\mathfrak g$ is therefore $\X(\mathbb R^2)$ and the
fundamental vector field $\ze_X(c)$ is now the vector field with value
$$ \ze_X(\ph) = \p_t| _0 (\ph \mapsto \ph \circ \on{exp}(tX)\circ \ph\i)
= \on{Ad}_\ph(X)$$
at the point $\ph$.
We now assume we have a positive definite inner product $G(X,Y)$ on the Lie
algebra $\X(\mathbb R^2)$ and that we use right translation to extend it to
a Riemannian metric on the full group $\on{Diff}(\mathbb R^2)$. This metric
being, by definition, invariant under the right group action, we have the
setting for momentum. The theory of the last section tells us to define the
momentum mapping by:
$$ j_X(\ph,Y) = G(\ze_X(\ph),Y).$$
Noether's theorem tells us that if $\ph(t)$ is a geodesic in
$\on{Diff}(\mathbb R^2)$ for this metric, then this momentum will be
constant along the lift of this geodesic to the tangent space. The lift of
$\ph(t)$, in the product decomposition of the tangent space is the curve:
$$ t \mapsto (\ph(t), \p_t(\ph) \circ \ph\i(t))$$
hence the theorem tells us that:
$$ G(\on{Ad}_{\ph(t)}(X), \p_t(\ph) \circ \ph\i(t)) = \on{constant}$$
for all $X$. If we further assume that Ad has an adjoint with respect to $G$:
$$ G(\on{Ad}_\ph(X),Y) \equiv G(X,\on{Ad}_\ph^*(Y))$$
then this invariance of momentum simplifies to:
$$ \boxed{\on{Ad}_{\ph(t)}^* \left(\p_t(\ph) \circ \ph\i(t)\right) = \on{constant}}$$
This is a very strong invariance and it encodes an integrated form of the
geodesic equations for the group.

\section{Geodesic equations and conserved momenta for almost local
Riemannian metrics}\label{almostlocalmetrics}

\subsection{The general almost local metric $G^\Ph$}
We have introduced above the $\Phi$-metrics:
$$ G^\Ph_c(h,k) 
:= \int_{S^1}\Ph(\ell_c,\ka_c(\th))\langle h(\th),k(\th) \rangle ds. $$
 
Since $\ell(c)$ is an integral operator the integrand is not a local
operator, but the nonlocality is very mild. We call it {\it almost local}.
The metric $G^\Ph$ is invariant under the reparametization group
$\on{Diff}(S^1)$ and under the Euclidean motion group.
Note (see \cite{MM1},~2.2) that 
\begin{align*}
D_{(c,h)}\ell_c &= \int_{S^1} \frac{\langle h_\th,c_\th \rangle}{|c_\th|}d\th 
=\int_{S^1} \langle D_s(h),v \rangle ds
\\&
=-\int_{S^1} \langle h, D_s(v) \rangle ds
=-\int_{S^1} \ka(c) \langle h, n\rangle ds
\\
D_{(c,h)}\ka_c &= \frac{\langle Jh_\th,c_{\th\th} \rangle}{|c_\th|^3} 
+ \frac{\langle Jc_\th,h_{\th\th} \rangle}{|c_\th|^3} 
-3\ka(c)\frac{\langle h_\th,c_\th \rangle}{|c_\th|^2} \\&
= \langle D_s^2(h), n \rangle -2 \ka \langle D_s(h), v \rangle.
\end{align*}
We compute the $G^\Ph$-gradients of $c\mapsto G^\Ph_c(h,k)$:
\begin{align*}
&D_{(c,m)}G_c^\Ph(h,k) = \int_{S^1}\Bigl(
 \p_1\Ph(\ell,\ka).D_{(c,m)}\ell_c.\langle  h,k\rangle 
 +\p_2\Ph(\ell,\ka).D_{(c,m)}\ka_c.\langle h,k \rangle
\\&\qquad\qquad\qquad\qquad\qquad\qquad\qquad\qquad
 +\Ph(\ell,\ka).\langle h,k\rangle \langle .D_s(m),v \rangle \Bigr)\,ds 
\\&= 
 -\int_{S^1} \ka_c \langle m,n\rangle ds\cdot
 \int_{S^1}\p_1\Ph(\ell,\ka)\langle  h,k\rangle ds
\\&\quad
+\int_{S^1}\Big(
\p_2\Ph(\ell,\ka) (\langle D_s^2(m),n \rangle 
  -2 \ka \langle D_s(m),v \rangle)
  +\Ph(\ell,\ka) \langle D_s(m),v\rangle \Bigr)\langle h, k \rangle ds 
\\&
=\int_{S^1}\Ph(\ell,\ka) \bigg\langle m\;,\; 
 \frac{1}{\Ph(\ell,\ka)}\bigg(
 -\ka_c\Big(\int\p_1\Ph(\ell,\ka)\langle  h,k\rangle ds\Big) n 
 + D_s^2\Big(\p_2\Ph(\ell,\ka) \langle h,k \rangle n \Big)
\\&\qquad\qquad\qquad\qquad
+ 2 D_s\Big(\p_2\Ph(\ell,\ka) \ka \langle h,k \rangle v \Big)
-D_s \Big( \Ph(\ell,\ka) \langle h,k \rangle v \Big) 
\bigg)\quad
\bigg\rangle ds
\end{align*}
According to \ref{setting} we should rewrite this as
\begin{align*}
D_{(c,m)}G^\Ph_c(h,k) &= 
G^\Ph_c(K^\Ph_c(m,h), k)
=G^\Ph_c\big(m, H^\Ph_c(h,k)\big),
\end{align*}
where the two $G^\Ph$-gradients $K^\Ph$ and $H^\Ph$ 
of $c\mapsto G^\Ph_c(h,k)$ are given by:
\begin{align*}
K^\Ph_c(m,h) &=
-\Big(\int_{S^1} \ka_c \langle m, n\rangle ds\Big)
 \frac{\p_1\Ph(\ell,\ka)}{\Ph(\ell,\ka)}h
\\&\quad
+\frac{\p_2\Ph(\ell,\ka)}{\Ph(\ell,\ka)}
  \Big(\langle D_s^2(m),n\rangle-2\ka\langle D_s(m),v\rangle\Big) h
+\langle D_s(m),v \rangle h
\\
H^\Ph_c(h,k) &=
\frac{1}{\Ph(\ell,\ka)}\bigg(
-\Big(\ka_c\int\!\p_1\Ph(\ell,\ka)\langle  h,k\rangle ds\Big) n 
+ D_s^2\Big(\p_2\Ph(\ell,\ka) \langle h,k \rangle n \Big)
+\\&\qquad\qquad\qquad\quad
+2D_s \Big(\p_2\Ph(\ell,\ka) \ka \langle h,k \rangle v \Big) 
- D_s\Big(\Ph(\ell,\ka)\langle h,k \rangle v \Big)\bigg)
\end{align*}

By substitution into the general formula of \ref{geodesic}, this gives the
geodesic equation for $G^\Ph$, but in a form which doesn't seem very
revealing, hence we omit it. Below we shall give the equation for the special
case of horizontal geodesics, i.e. geodesics in $B_i$.

\subsection{Conserved momenta for $G^\Ph$}\label{momentumPh}
According to \ref{momentum} the momentum mappings 
for the reparametrization, translation and rotation group actions 
are conserved along any geodesic $t\mapsto c(t,\;.\;)$:
$$ \boxed{
\begin{aligned}
&\Ph(\ell_c,\ka_c)\langle v,c_t \rangle |c_\th|^2\;\in\X(S^1)
&\quad&\text{reparametrization momentum}
\\
&\int_{S^1}\Ph(\ell_c,\ka_c) c_t ds \;\in \mathbb R^2
&\quad&\text{linear momentum}
\\
&\int_{S^1}\Ph(\ell_c,\ka_c) \langle Jc,c_t \rangle ds\in \mathbb R
  &\quad&\text{angular momentum}
\end{aligned}} $$

Note that setting the reparametrization momentum to 0 and doing symplectic
reduction there amounts exactly to investigating the quotient space
$B_i(S^1,\mathbb R^2)=\on{Imm}(S^1,\mathbb R^2)/\on{Diff}(S^1)$ and using
horizontal geodesics for doing so; a horizontal geodesic is one 
for which $\langle v,c_t \rangle=0$; or equivalently it
is $G^\Ph$-normal to the $\on{Diff}(S^1)$-orbits. If it is normal at one time
it is normal forever (since the reparametrization momentum is conserved).
This was the approach taken in \cite{MM1}.

\subsection{Horizontality for $G^\Ph$}
The tangent vectors to the $\on{Diff}(S^1)$ orbit through $c$ are 
$T_c(c\o\on{Diff}(S^1)) = \{X.c_\th:X\in C^\infty(S^1,\mathbb R)\}$. 
Thus the bundle of horizontal
vectors is
\begin{align*}
\mathcal N_c & =\{h\in C^\infty(S^1,\mathbb R^2): \langle h,v\rangle =0\}
\\
&=\{a.n\in C^\infty(S^1,\mathbb R^2): a\in C^\infty(S^1,\mathbb R)\}
\end{align*}
A tangent vector $h\in T_c\on{Imm}(S^1,\mathbb R^2)=C^\infty(S^1,\mathbb R^2)$ 
has an orthonormal de\-compo\-si\-tion 
\begin{align*}
h&=h^\top+h^\bot \in T_c(c\o\on{Diff}^+(S^1)) \oplus \mathcal N_c\qquad
  \text{ where}
\\
h^\top &= \langle h,v \rangle v 
  \in  T_c(c\o\on{Diff}^+(S^1)),
\\
h^\bot &= \langle h, n \rangle n \in  \mathcal N_c,
\end{align*}
into smooth tangential and normal components, independent of the choice of
$\Ph(\ell,\ka)$.
For the following result the proof given in \cite{MM1},~2.5 works without
any change:

\begin{lem*}
For any smooth path $c$ in $\on{Imm}(S^1,\mathbb R^2)$ there exists a smooth
path $\ph$ in $\on{Diff}(S^1)$ with $\ph(0,\;.\;)=\on{Id}_{S^1}$ depending
smoothly on $c$ such that
the path $e$ given by $e(t,\th)=c(t,\ph(t,\th))$ is horizontal:
$e_t\bot e_\th$. \qed
\end{lem*}

Consider a path $t\mapsto c(\cdot,t)$ in the manifold $\on{Imm}(S^1,\mathbb
R^2)$. It projects to a path $\pi \circ c$ in $B_i(S^1,\mathbb R^2)$ whose
energy is
called the {\it horizontal energy} of $c$:
\begin{gather*}
{\begin{aligned}
E_{G^\Ph}^{\text{hor}}(c) &= 
E_{G^\Ph}(\pi\o c) 
= \tfrac12\int_a^b G^\Ph_{\pi(c)}(T_c\pi.c_t,T_c\pi.c_t)\,dt
\\& 
= \tfrac12\int_a^b G^\Ph_{c}(c_t^\bot,c_t^\bot)\,dt
=\tfrac12\int_a^b\int_{S^1}\Ph(\ell_c,\ka_c) 
  \langle c_t^\bot,c_t^\bot\rangle ds\,dt
\end{aligned}}
\\
\boxed{\quad
E^{\text{hor}}_{G^\Ph}(c) 
=\tfrac12\int_a^b\int_{S^1} \Ph(\ell_c,\ka_c)
\langle c_t,n \rangle^2 \,d\th dt
\quad}
\end{gather*}
For a horizontal path this is just the usual energy.
As in \cite{MM1},~3.12
we can express $E^{\text{hor}}(c)$ as an integral over the
graph $S$ of $c$, the immersed surface
$S\subset \mathbb R^3$ parameterized by $(t,\th)\mapsto
(t,c(t,\th))$, in terms of the surface area
$d\mu_S = |\Ph_t\x \Ph_\th|\,d\th\,dt$
and the unit normal $n_S=(n_S^0,n_S^1,n_S^2)$ of $S$:
\begin{equation*}
E^{\text{hor}}_{G^\Ph}(c) 
=\tfrac12\int_{[a,b]\x S^1}\Ph(\ell_c,\ka_c)
  \frac{|n_S^0|^2}{\sqrt{1-|n_S^0|^2}}\,d\mu_S
\end{equation*}
Here the final expression is only in terms of the surface $S$ and its
fibration over the time axis, and is valid for any path $c$.
This anisotropic area functional has to be minimized in order to prove that 
geodesics exists between arbitrary curves (of the same degree) in
$B_{i}(S^1,\mathbb R^2)$.

\subsection{The horizontal geodesic equation}
\label{horizgeod}
Let $c(\th,t)$ be a horizontal geodesic for the metric $G^\Ph$. Then
$c_t(\th,t) = a(\th,t).n(\th,t)$. 
Denote the integral of a function over the curve with
respect to arclength by a bar. Then the geodesic equation for horizontal
geodesics is:
$$ \boxed{
\begin{aligned} 
a_t = \frac{-1}{2\Ph} \Big(&\left( -\ka \Ph 
+ \ka^2 \p_2\Ph \right)a^2 -D_s^2\left(\p_2\Ph \cdot a^2\right) 
+ 2\p_2\Ph \cdot aD_s^2(a) 
\\& -2\p_1\Ph \cdot\overline{(\ka a)}\cdot a 
+ \overline{(\p_1\Ph \cdot a^2)}\cdot \ka \Big)
\end{aligned}}
$$

This comes immediately from the formulas for $H$ and $K$ in the metric
$G^\Ph$ when you substitute $m=h=k=a.n$ and consider only the $n$-part. We
obtain in this case:
\begin{align*}
\Ph \cdot \langle K,n \rangle 
&= -\overline{(\ka a)}. \p_1\Ph .a + \p_2\Ph .D_s^2(a).a 
+ \p_2\Ph .\ka^2 a^2 - \Ph \ka a^2\\
\Ph \cdot \langle H,n \rangle 
&= -\overline{(\p_1\Ph a^2)}. \ka + D_s^2(\p_2\Ph . a^2) 
+\p_2\Ph. \ka^2 a^2 - \Ph \ka a^2.
\end{align*}
and the geodesic formula follows by substitution.

\subsection{Curvature on $B_{i,f}(S^1,\mathbb R^2)$ for $G^\Ph$}
\label{curvatureG^Ph}
We compute the curvature of $B_i(S^1,\mathbb R^2)$ in the general almost
local metric $G^\Ph$. We proceed as in \cite{MM1},~{2.4.3}.
We use the following chart near $C\in B_{i}(S^1,\mathbb R^2)$.
Let $c\in\on{Imm}_{f}(S^1,\mathbb R^2)$ be parametrized by arclength with 
$\pi(c)=C$ of length $L$, with unit normal $n_c$. 
We assume that the parameter $\th$ runs in the scaled circle $S^1_L$ below.
\begin{align*}
&\ps: C^\infty(S^1_L,(-\ep,\ep)) \to \on{Imm}_{f}(S^1_L,\mathbb R^2),
  \qquad \mathcal Q(c):= \ps(C^\infty(S^1_L,(-\ep,\ep)))\\
&\ps(f)(\th) = c(\th) + f(\th)n_c(\th) = c(\th) + f(\th)ic'(\th),\\
&\pi\o\ps :C^\infty(S^1_L,(-\ep,\ep)) \to B_{i,f}(S^1,\mathbb R^2),
\end{align*}
where $\ep$ is so small that $\ps(f)$ is an embedding for each $f$.
We have (see \cite{MM1}, {2.4.3})
\begin{align*}
\ps(f)' &= c' + f'ic' + fic'' = (1-f\ka_c)c' + f'ic' \\
\ps(f)'' &= c'' + f''ic' + 2f'ic'' + fic'''
     = -(2f'\ka_c+f\ka_c')c'+ (\ka_c+f''-f\ka_c^2)ic'  \\
n_{\ps(f)} &= 
     \frac1{\sqrt{(1-f\ka_c)^2+{f'}^2}}\Bigl((1-f\ka_c)ic'-f'c'\Bigr),\\
T_f\ps.h &= h.ic'\qquad \in C^\infty(S^1,\mathbb R^2) = 
     T_{\ps(f)}\on{Imm}_{f}(S^1_L,\mathbb R^2)\\
&=\frac{h(1-f\ka_c)}{\sqrt{(1-f\ka_c)^2+{f'}^2}}n_{\ps(f)}
     +\frac{hf'}{(1-f\ka_c)^2+{f'}^2}\ps(f)',\\
(T_f\ps.h)^\bot &=\frac{h(1-f\ka_c)}{\sqrt{(1-f\ka_c)^2+{f'}^2}}n_{\ps(f)}
     \quad\in \mathcal N_{\ps(f)},\\
\ka_{\ps(f)}&= \frac{1}{((1-f\ka_c)^2+{f'}^2)^{3/2}}
     \langle i\ps(f)',\ps(f)'' \rangle \\
&= \ka_c+(f''+f\ka_c^2)+(f^2\ka_c^3+\tfrac12{f'}^2\ka_c+ff'\ka_c'+2ff''\ka_c)
     +O(f^3)
\\
\ell(\ps(f)) &= \int_{S^1_L}|\ps(f)|\,d\th =
  \int_{S^1_L}(1-2f\ka_c + f^2\ka_c^2+{f'}^2)^{1/2}\,d\th
\\&
= \int_{S^1_L}\Big(1-f\ka_c +\frac{{f'}^2}{2} +O(f^3)\Big)\,d\th
= L- \overline{f\ka_c} +\tfrac12\overline{{f'}^2} +O(f^3)
\end{align*}
where we use the shorthand
$\overline{g}=\int_{S^1_L}g(\th)\,d\th=\int_{S^1_L}g(\th)\,ds$.
Let $G^\Ph$ denote also the induced metric on 
$B_{i,f}(S^1_L,\mathbb R^2)$. 
Since $\pi$ is a Riemannian submersion, 
for 
$f\in C^\infty(S^1_L,(-\ep,\ep))$ and $h,k\in C^\infty(S^1_L,\mathbb R)$ we
have
\begin{align*}
&((\pi\o\ps)^*G^\Ph)_f(h,k) = 
G^\Ph_{\pi(\ps(f))}\Bigl(T_f(\pi\o\ps)h,T_f(\pi\o\ps)k\Bigr)\\
&=G^\Ph_{\ps(f)}\Bigl((T_f\ps.h)^\bot,(T_f\ps.k)^\bot\Bigr)
=\int_{S^1_L}\Ph(\ell(\ps(f)),\ka_{\ps(f)})
     \frac{hk(1-f\ka_c)^2}{\sqrt{(1-f\ka_c)^2+{f'}^2}}\;d\th\\
\end{align*}
We have to compute second derivatives in $f$ of this.
For that we expand the main contributing 
expressions in $f$ to order 2:
\begin{align*}
\\
&(1-f\ka)^2(1-2f\ka+f^2\ka^2+{f'}^2)^{-1/2} 
  =1-f\ka-\tfrac12{f'}^2+O(f^3)
\\&
\Ph(\ell,\ka) =
\Ph(L,\ka_c) + \p_1\Ph(L,\ka_c)(\ell-L) 
+ \p_2\Ph(L,\ka_c)(\ka-\ka_c) 
\\&\quad
+\p_1\p_2\Ph(L,\ka_c)(\ell-L)(\ka-\ka_c) 
\\&\quad
+\frac{\p_1^2\Ph(L,\ka_c)}{2}(\ell-L)^2 +
\frac{\p_2^2\Ph(L,\ka_c)}{2}(\ka-\ka_c)^2 + O(3)
\end{align*}
We simplify notation as $\ka=\ka_c$, $\Ph=\Ph(L,\ka_c)$,
$((\pi\o\ps)^*G^\Ph)_f=G^\Ph_f$ etc.\ and expand the metric: 
\begin{align*}
G^\Ph_f(h,k)
=\int_{S^1_L}hk\bigg(&
\Ph - \p_1\Ph.\overline{f\ka} 
+\p_2\Ph.(f''+f\ka^2)
-\Ph.f\ka 
\\&
+\tfrac12\p_1\Ph.\overline{{f'}^2}
+\p_2\Ph.(f^2\ka^3+\tfrac12{f'}^2\ka+ff'\ka'+2ff''\ka)
\\&
-\p_1\p_2\Ph.\overline{f\ka}(f''+f\ka^2)
+\frac{\p_1^2\Ph}{2}(\overline{f\ka})^2 
+\frac{\p_2^2\Ph}{2}(f''+f\ka^2)^2
\\&
+ \p_1\Ph.f\ka.\overline{f\ka}
-\p_2\Ph.f\ka.(f''+f\ka^2)
-\Ph.\tfrac12{f'}^2 
\bigg)d\th      + O(f^3)
\end{align*}
Note that $G^\ph_0(h,k)=\int_{S^1_L}hk\Ph\,d\th$.
We differentiate the metric
and compute the Christoffel symbol at the center $f=0$
\begin{align*}
-2G^A_0(&\Ga_0(h,k),l) = -dG^A(0)(l)(h,k) +dG^A(0)(h)(k,l) +dG^A(0)(k)(l,h) 
\\
=\int_{S^1_L}\Big(&
-\p_1\Ph.\overline{h\ka}.kl-\p_1\Ph.h.\overline{k\ka}.l
+\p_1\Ph.hk\int l\ka\,d\th_1
-\p_2\Ph''.hkl
\\&
-2\p_2\Ph'.h'kl-2\p_2\Ph'.hk'l-2\p_2\Ph.h'k'l
+\p_2\Ph.hkl\ka^2 -\Ph.hkl\ka 
\Big)d\th   
\end{align*}
Thus 
\begin{align*}
\Ga_0(h,k)
&= \frac1{2\Ph}\Big(\p_1\Ph.(\overline{h\ka}.k+h.\overline{k\ka})
-\ka\overline{\p_1\Ph.hk}
\\&\quad
+\p_2\Ph''.hk+2\p_2\Ph'.h'k+2\p_2\Ph'hk'
+2\p_2\Ph.h'k'
\\&\quad
-\p_2\Ph.hk\ka^2 +\Ph.hk\ka\Big) 
\end{align*}
Letting $h=k=f_t=a$, this leads to the geodesic equation
from \ref{horizgeod}.
For the sectional curvature we use the following formula which is valid in a chart:
\begin{align*}
&2R_f(m,h,m,h)=2G^A_f(R_f(m,h)m,h) =
\\&
= -2d^2G^A(f)(m,h)(h,m)  +d^2G^A(f)(m,m)(h,h)  +d^2G^A(f)(h,h)(m,m)  
\\&\quad
-2G^A(\Ga(h,m),\Ga(m,h))
+2G^A(\Ga(m,m),\Ga(h,h))
\end{align*}
The sectional curvature at the two-dimensional subspace $P_f(m,h)$ of the
tangent space which is spanned by $m$ and $h$ is then given by:
\begin{equation*}
k_f(P_f(m,h)) = - \frac{G^\Ph_f(R(m,h)m,h)}{\|m\|^2\|h\|^2-G^\Ph_f(m,h)^2}.
\end{equation*}
We compute this directly for $f=0$,
using the expansion up to order 2 of $G^A_f(h,k)$ and the Christoffels.
We let $W(\th_1,\th_2) = h(\th_1)m(\th_2)-h(\th_2)m(\th_1)$ so that its
second derivative 
$\p_2W(\th_1,\th_1)=W_2(\th_1,\th_1)=h(\th_1)m'(\th_1)-h'(\th_1)m(\th_1)$ 
is the Wronskian of $h$ and $m$.
Then we have
our final result for the main expression in the horizontal sectional
curvature, where we use $\int=\int_{S^1_L}$, $\overline{g}=\int_{S^1_L}
g\,ds$, and $\Ph_1=\p_1\Ph$ etc. Also recall that the base curve is
parametrized by arc-length.
$$\boxed{
\begin{aligned}
&R^\Ph_0(m,h,m,h)=G^\Ph_0(R_0(m,h)m,h) =
\\&
= \int \Big(\ka.\Ph_2-\frac{\Ph}{2} +\frac{\Ph_2.\Ph_2''
  -2(\Ph_2')^2-(\Ph_2\ka)^2}{2\Ph}\Big)(\th_1)W_2(\th_1,\th_1)^2\,d\th_1
\\&
+ \int \frac{\Ph_{22}(\th_1)}{2}W_{22}(\th_1,\th_1)^2\,d\th_1
\\&
+ \int \Big(\frac{\Ph_1'\Ph_2}{\Ph}-\frac{\Ph_1\Ph_2\Ph_1'}{\Ph^2}\Big)(\th_1)
       W_2(\th_1,\th_1)\int W(\th_1,\th_2)\ka(\th_2)\,d\th_2\,d\th_1
\\&
+ \int \Big(\frac{\Ph_1\Ph_2}{\Ph}-\Ph_{12}\Big)(\th_1)
       W_{22}(\th_1,\th_1)\int W(\th_1,\th_2)\ka(\th_2)\,d\th_2\,d\th_1
\\&
+ \iint \frac{\Ph_1(\th_1)}{2}\Big(1-\frac{\Ph_2.\ka}{\Ph}(\th_2)\Big)
  W_1(\th_1,\th_2)^2\,d\th_2\,d\th_1
\\&
+ \iint \Big(\frac{\Ph_2.\ka^3-\Ph_2''.\ka}{4\Ph} -\frac{\ka^2}{4}
  +\Big(\frac{\Ph_2'.\ka}{2\Ph}\Big)'
  +\overline{\big(\frac{\ka^2}{8\Ph}\big)}.\Ph_1\Big)(\th_1)
  \Ph_1(\th_2)W(\th_1,\th_2)^2\,d\th_2\,d\th_1
\\&
+ \iiint \Big(\frac{\Ph_{11}}{2}-\frac{\Ph_1^2}{4\Ph}\big)(\th_1)
  -\Ph_1(\th_1)\frac{\Ph_1}{2\Ph}(\th_2)\Big)
\\&\qquad\qquad\qquad\qquad
  \ka(\th_2)\ka(\th_3)W(\th_1,\th_2)W(\th_1,\th_3)\,d\th_2\,d\th_1\,d\th_3
\end{aligned}
}$$

\subsection{Special case: the metric $G^A$}
If we choose $\Ph(\ell_c,\ka_c)=1+A\ka_c^2$ then we obtain 
the metric used in \cite{MM1}, given by  
$$
G^A_c(h,k) = \int_{S^1}(1+A\ka_c(\th)^2)\langle h(\th),k(\th) \rangle ds.
$$
As shown in our earlier paper, $\sqrt{\ell}$ is Lipschitz in this metric and the metric dominates the Frechet metric.

The horizontal geodesic equation for the $G^A$-metric reduces to
$$ a_t = \frac{-\tfrac12 \ka_c a^2 + A \left( 
a^2(-D_s^2(\ka_c) + \tfrac12 \ka_c^3) - 4 D_s(\ka_c) a D_s(a) 
- 2\ka_c D_s(a)^2 \right)}{1+A\ka_c^2} $$
as found in \cite{MM1},4.2. Along a geodesic
$t\mapsto c(t,\;.\;)$ we have the following conserved quantities:
\begin{align*}
&(1+A\ka_c^2)\langle v,c_t \rangle |c_\th|^2 \;\in\X(S^1)
&\quad&\text{reparametrization momentum} \\
&\int_{S^1}(1+A\ka_c^2) c_t ds \;\in \mathbb R^2
&\quad&\text{linear momentum} \\
&\int_{S^1}(1+A\ka_c^2) \langle Jc,c_t \rangle ds \; \in \mathbb R
  &\quad&\text{angular momentum}
\end{align*}
For $\Ph(\ell,\ka)=1+A\ka^2$ we have $\p_1\Ph=0$, $\p_2\Ph=2A\ka$,
$\p_2^2\Ph=2A$, and the general curvature formula in \ref{curvatureG^Ph}
for the horizontal curvature specializes to the formula in \cite{MM1}, 4.6.4:
$$ R^\Ph_0(m,h,m,h)
= \int \Big(-\frac{(1-A\ka^2)^2 -4A^2 \ka \ka'' + 8A^2\ka'^2}{2(1+A\ka^2)}
W_2^2 + A W_{22}^2\Big) d\th. $$

\subsection{Special case: the conformal metrics}
We put $\Ph(\ell(c),\ka(c))=\Ph(\ell(c))$ and obtain the
metric proposed by Menucci and Yezzi and, for $\Phi$ linear, independently 
by Shah \cite{Shah}:
$$
G^\Ph_c(h,k)= \Ph(\ell_c)\int_{S^1} \langle  h,k\rangle ds =
\Ph(\ell_c)G^0_c(h,k).
$$
All these metrics are conformally equivalent to the basic $L^2$-metric
$G^0$.
As they show, the infimum of path lengths in this metric is positive so
long as $\Ph$ satifies an inequality $\Ph(\ell) \ge C.\ell$ for some $C>0$.
This follows, as in \cite{MM1}, 3.4, by the inequality on area swept out by
the curves in a horizontal path $c_t = a.n$:
\begin{align*}
 \int |a|.ds &\le \left( \int a^2.ds \right)^{1/2} \cdot \ell^{1/2} \le
 \left(\frac{\ell}{\Ph(\ell)} \right)^{1/2}\cdot (G^\Ph(a,a))^{1/2} \\
\text{Area swept out} &\le
\max_t\left(\frac{\ell_{c(t,\cdot)}}{\Ph(\ell_{c(t,\cdot)}}
\right)^{1/2}\cdot \left( G^\Ph\text{-path length} \right) \le
\frac{G^\Ph\text{-path length}}{\sqrt{C}}.
\end{align*}
The horizontal geodesic equation reduces to:
$$ a_t = \frac{\ka}{2} a^2 - \frac{\p_1\Ph}{\Ph} \cdot \bigg( \tfrac12
\left(\int a^2.ds\right) \ka -\left( \int \ka . a . ds \right) a \bigg)
$$
If we change variables and write $b(s,t) = \Ph(\ell(t)).a(s,t)$, then this
equation simplifies to:
$$ \begin{boxed}
{\quad b_t = \frac{\ka}{2\Ph} \left( b^2 - \frac{\p_1 \Ph}{\Ph} \int b^2
\right) \quad}
\end{boxed}$$

Along a geodesic $t\mapsto c(t,\;.\;)$ we have the following conserved
quantities:
\begin{align*}
&\Ph(\ell_c)\langle v,c_t \rangle |c'(\th)|^2\;\in\X(S^1)
&\quad&\text{reparametrization momentum}
\\
&\Ph(\ell_c)\int_{S^1} c_t ds\;\in \mathbb R^2
&\quad&\text{linear momentum}
\\
&\Ph(\ell_c)\int_{S^1}\langle Jc,c_t \rangle ds \in \mathbb R
  &\quad&\text{angular momentum}
\end{align*}

For the conformal metrics, sectional curvature has been computed by Shah
\cite{Shah} using the method of local charts from \cite{MM1}. 
We specialize formula \ref{curvatureG^Ph} to the case that 
$\Ph(\ell,\ka)=\Ph(\ell)$ is independent of $\ka$. Then $\p_2\Ph=0$.
We also assume that $h,m$ are orthonormal so that $\Ph\overline{h^2}=\Ph\overline{m^2}=1$
and $\Ph\overline{hm}=0$. Then the 
the sectional curvature at the two-dimensional subspace $P_0(m,h)$ of the
tangent space which is spanned by $m$ and $h$ is then given by:
\begin{align*}
&k_0(P_0(m,h)) = - \frac{G^\Ph_0(R_0(m,h)m,h)}{\|m\|^2\|h\|^2-G^\Ph_0(m,h)^2}=
\\&=
\tfrac12\Ph.\overline{W(h,m)^2}
+\frac{\p_1\Ph}{4\Ph}.(\overline{m^2\ka^2}
  +\overline{h^2\ka^2})
+ \frac{3(\p_1\Ph)^2-2\Ph\p_1^2\Ph}{4\Ph^2}(\overline{h\ka}^2+\overline{m\ka}^2)
\\&\quad
-\frac{\p_1\Ph}{2\Ph}(\overline{{m'}^2}+\overline{{h'}^2})
-\frac{(\p_1\Ph)^2}{4\Ph^3}\overline{\ka^2}
\end{align*}
which is the same as the equation (11) in \cite{Shah}. 
Note that the first line is positive while the last line is negative.
The first term is the curvature term for the $H^0$-metric.
The key point about this formula is how many positive terms it has. This
makes it very hard to get smooth geodesics in this metric. For example, in
the case where $\Ph(\ell)=c.\ell$, the analysis of Shah \cite{Shah} proves that
the infimum of $G^\Ph$ path length between two embedded curves $C$ and $D$
is exactly the area of the symmetric difference of their interiors:
$\on{Area}(\on{Int}(C)\De \on{Int}(D))$, but that this length is realized
by a smooth path if and only if $C$ and $D$ can be connected by
`grassfire', i.e.\ a family in which the length $|c_t(\th,t)| \equiv 1$.

\subsection{Special case: the smooth scale invariant metric $G^{SI}$}
\label{scalemetric}
Choosing the function $\Ph(\ell,\ka)=\ell^{-3}+A\frac{\ka^2}{\ell}$ we
obtain the metric:
$$ G^{SI}_c(h,k)
  = \int_{S^1} \Big(\frac1{\ell_c^3}+A\frac{\ka_c^2}{\ell_c}\Big)
  \langle h,k\rangle ds. $$
The beauty of this metric is that (a) it is scale invariant and (b)
$\log(\ell)$ is Lipschitz, hence the infimum of path lengths is always
positive. Scale invariance is clear: changing $c,h,k$ to $\la\cdot c, \la
\cdot h, \la\cdot k$ changes $\ell$ to $\la\cdot \ell$ and $\ka$ to
$\ka/\la$ so the $\la$'s in $G^{SI}$ cancel out. To see the second fact,
take a horizontal path $c_t = a\cdot n,\; 0 \le t \le 1$, and abbreviate
the lengths of the curves in this path, $\ell_{c(t,\cdot)}$, to $\ell(t)$.
Then we have:
\begin{align*}
\frac{\partial \log \ell(t)}{\partial t} &= \frac{1}{\ell(t)} \int_{S^1}
\ka_{c(t,\cdot)}(\th) \cdot a(\th,t) ds, \quad \text{hence} \\
\left| \frac{\partial \log \ell(t)}{\partial t}\right| &= \left( \frac{\int
\ka^2 a^2 ds}{\ell(t)} \right)^{1/2} \cdot \left( \frac{ \int 1\cdot
ds}{\ell(t)} \right)^{1/2}\\
&\le \frac{1}{\sqrt{A}} \left(G^{SI}(a,a)\right)^{1/2}, \quad
\text{hence}\\
|\log(\ell(1))-\log(\ell(0))| &\le \text{SI-path length}/\sqrt{A}.
\end{align*}
Thus in a path whose length in this metric is $K$, the lengths of the
individual curves can increase or decrease at most by a factor
$e^{K/\sqrt{A}}$. Now use the same argument as above to control the area
swept out by such a path:
\begin{align*}
\int |a| ds & \le \left( \int a^2 ds \right)^{1/2} \cdot \left( \int 1
\cdot ds \right)^{1/2} \\
& \le \left( \ell^3 G^{SI}(a,a)\right)^{1/2} \cdot \ell^{1/2} = \ell^2
\cdot G^{SI}(a,a)^{1/2}, \quad \text{hence}\\
\text{Area-swept-out} &\le e^{K/\sqrt{A}} \ell(0)^2 \cdot K
\end{align*}
which verifies the second fact.
We can readily calculate the geodesic equation for horizontal geodesics in
this metric as another special case of the equation for $G^\Ph$:
\begin{align*}
a_t &= \frac{-1}{1+A(\ell\ka)^2}\bigg( \left( -1+A(\ell\ka)^2\right)
\frac{\ka a^2}{2}-A\ell^2D_s^2(\ka)a^2 -2A\ell^2\ka D_s(a)^2 
\\&\quad 
- 4A\ell^2 D_s(\ka) a D_s(a) 
+  \left( 3+A(\ell\ka)^2\right) \overline{(a \ka)}  \cdot a -
\frac{3}{2} \overline{(a^2)} \cdot \ka - \frac{A\ell^2}{2} \overline{(\ka
a)^2} \cdot \ka \bigg)
\end{align*}
where the ``overline'' stands now for the {\it average} of a function over
the curve, i.e. $\int \cdots ds/\ell$.
Since this metric is scale invariant, there are now {\it four} conserved
quantities, instead of three:
\begin{align*}
&\Ph(\ell, \ka)\langle v,c_t \rangle |c'(\th)|^2\;\in\X(S^1)
&\quad&\text{reparametrization momentum}\\
&\int_{S^1}\Ph(\ell, \ka) c_t ds\;\in \mathbb R^2
&\quad&\text{linear momentum}\\
&\int_{S^1}\Ph(\ell, \ka)\langle Jc,c_t \rangle ds \in \mathbb R
  &\quad&\text{angular momentum}\\
&\int_{S^1}\Ph(\ell, \ka)\langle c,c_t \rangle ds \in \mathbb R
  &\quad&\text{scaling momentum}
\end{align*}
It would be very interesting to compute and compare geodesics in these
special metrics.

\subsection{The Wasserstein metric and a related $G^\Ph$-metric}
The Wasserstein metric (also known as the Monge-Kantorovich metric) is a
metric between probability measures on a common metric space, see
\cite{ambrosio}, and \cite{ambrosio-book} for more details. It has been
studied for many years globally and is defined as follows: let $\mu$ and
$\nu$ be 2 probability measures on a metric space $(X,d)$. Consider all
measures $\rh$ on $X \times X$ whose marginals under the 2 projections are
$\mu$ and $\nu$. Then:
$$ d_{\on{wass}}(\mu,\nu) = \inf_{\rh: p_{1,*}(\rh)=\mu, p_{2,*}(\rh)=\nu}
\iint_{X \times X} d(x,y) d\rh(x,y).$$
It was discovered only recently by Benamou and Brenier \cite{BenBren} that, 
if $X = \mathbb R^n$, this is, in fact, path length for a Riemannian metric 
on the space of probability measures $\mathcal P$. In their theory, the tangent
space at $\mu$ to the space of probability measures and the infinitesimal
metric are defined by:
$$ T_{\mu,\mathcal P} = \left\{\text{vector fields } h = \nabla f \text{
completed in the norm } \int |h |^2 d\mu \right\}$$
where the tangent $h$ to a family $t\mapsto \mu(t)$ is defined by the identity:
$$ \frac{\p \mu}{\p t} + \on{div} (h . \mu) = 0.$$
In our case, we want to assign to an immersion $c$ the scaled arc length
measure $\mu_c = ds / \ell$. This maps $B_i$ to $\mathcal P$. The claim is
that the pull-back of the Wasserstein metric by this map is intermediate
between $G^{\ell\i}$ and $G^{\Ph_W}$, where
$$\Ph_W(\ell,\ka)= \ell\i+\tfrac{1}{12}\ell\ka^2.$$
This is not hard to work out.
\begin{enumerate}
\item Because we are mod-ing out by vector fields of norm 0, the vector
field $h$ is defined only along the curve $c$ and its norm is $\ell\i .\int
\|h\|^2 ds$.
\item If we split $h = av+bn$, then the condition that $h=\nabla f$ means
that $\int a.ds = 0$ and the norm is $\ell\i.\int (a^2+b^2)ds$.
\item But moving c infinitesimally by $h$, scaled arc length
parametrization of $c$ must still be scaled arc length. Let $c(\cdot,t) =
c+th$. Then this means $|c_\th|_t = \on{cnst.}|c_\th|$ at $t=0$. Since
$|c_\th|_t = \langle c_{t\th}, c_\th \rangle / |c_\th|$, this condition is
the same as $\langle D_s(av+bn),v \rangle = \on{cnst.}$, or $D_s a-b\ka_c =
\on{cnst.}$.
\item Combining the last 2 conditions on $b$, we get a formula for $a$ in
terms of $b$, namely $a = K *( b\ka_c)$, where we convolve with respect to
arc length using the kernel $K(x) = \on{sign}(x)/2 - x/\ell, -\ell \le x
\le \ell$.
\item Finally, since $|K*f|(x) \le |K|.|f| = \sqrt{\ell/12}|f|$ for all
$f$, it follows that
$$\ell\i.\int b^2ds \le \ell\i.\int (a^2+b^2)ds \le \ell\i.\int (b^2 +
\frac{(\ell\ka)^2}{12}.b^2)ds$$
which sandwiches the Wasserstein norm between $G^{\ell\i}$ and $G^{\Ph_W}$
for $\Ph_W= \ell\i.(1+(\ell\ka)^2/12)$.
\end{enumerate}

\section{Immersion-Sobolev metrics on $\on{Imm}(S^1,\mathbb R^2)$ and on
$B_i$}
\label{G^{imm,n}}

\subsection{The $G^{\on{imm},n}$-metric}
\label{H^n-geodesic}
We note first that the differential operator $D_s=\frac{\p_\th}{|c_\th|}$
is anti self-adjoint for the metric $G^0$, i.e., for all $h,k\in
C^\infty(S^1.\mathbb R^2)$ we have
$$ \int_{S^1}\Big\langle D_s(h),k \Big\rangle ds
=\int_{S^1}\Big\langle h,-D_s(k) \Big\rangle ds $$
We can define a Sobolev-type weak Riemannian metric\footnote{There are
other choices for the higher order terms, e.g.\ summing all the
intermediate derivatives with or without binomial coefficients. These
metrics are all equivalent and the one we use leads to the simplest
equiations.} on $\on{Imm}(S^1,\mathbb R^2)$ which is invariant under the
action of $\on{Diff}(S^1)$ by:
\begin{align*}
G^{\on{imm},n}_c(h,k) &=
\int_{S^1} \left(\langle h,k \rangle + A.\langle D_s^n h,
   D_s^n k  \rangle \right). ds
\tag{1}\\&
= \int_{S^1}\langle L_n(h),k\rangle ds
  \qquad\text{  where }
\\
L_n(h) \text{ or } L_{n,c}(h) &= (I + (-1)^n A.D_s^{2n})(h)
\tag{2}\end{align*}

The interesting special case $n=1$ and $A \rightarrow \infty$ has been
studied by Trouv\'e and Younes in \cite{YT, Y} 
and by Mio, Srivastava and Joshi in \cite{MS, MS2}. In this case, the
metric reduces to:
$$ G^{\on{imm},1,\infty}_c(h,k) = \int_{S^1} \langle D_s(h), D_s(k)
\rangle.ds$$
which ignores translations, i.e.\ it is a metric on $\on{Imm}(S^1,\mathbb
R^2) \text{ modulo translations}$. Now identify $\mathbb R^2$ with $\mathbb
C$, so that this space embeds as follows:
\begin{align*} 
\on{Imm}(S^1,\mathbb R^2)/\rm{transl.} &\hookrightarrow C^\infty(S^1, \mathbb C) \\
c & \longmapsto c_\th.
\end{align*}
Then Trouv\'e and Younes use the new shape space coordinates $Z(\th) =
\sqrt{c_\th(\th)}$ and Mio et al use the coordinates $\Ph(\th) =
\log(c_\th(\th))$ -- with {\it complex} square roots and logs. Both of
these unfortunately require the introduction of a discontinuity, but this
will drop out when you minimize path length with respect to
reparametrizations. The wonderful fact about $Z(\th)$ is that in a family
$Z(t,\th)$, we find:
$$ 
Z_t =\frac{c_{t,\th}}{2\sqrt{c_\th}}, \quad \text{so } 
\int_{S^1} |Z_t|^2 d\th = \tfrac14 \int \frac{|c_{t,\th}|^2}{|c_\th|^2} |c_\th| d\th 
= \tfrac14 \int |D_s(c_t)|^2 ds
$$
so the metric becomes a {\it constant} metric on the vector space of
functions $Z$. With $\Ph$, one has $\int |\Ph_t|^2 ds = \int
|D_s(c_t)|^2ds$, which is simple but not quite so nice. One can expect a
very explicit representation of the space of curves in this metric.

Returning to the general case, for each fixed $c$ of length $\ell$, the
differential operator $L_{n,c}$ is simply the constant coefficient ordinary
differential operator $f \mapsto f +(-1)^n A .f^{(2n)}$ on the $s$-line
modulo $\ell.\mathbb Z$. Thus its Green's function is a linear combination
of the exponentials $\exp(\la.x)$, where $\la$ are the roots of $1+(-1)^n
A.\la^{2n} = 0$. A simple verification gives its Green's function (which we
will not use below):
$$ F_n(x) = \frac{1}{2n}\cdot \sum_{\la^{2n} = (-1)^{n+1}/A}
\frac{\la}{1-e^{\la \ell}} e^{\la x}, \quad 0 \le x \le \ell. $$
This means that the dual metric 
$\check G^{\on{imm},n}_c=(G^{\on{imm},n}_c)\i$ 
on the {\it smooth cotangent space}
$C^\infty(S^1,\mathbb R^2)\cong G^0_c(T_c\on{Imm}(S^1,\mathbb R^2))
\subset T_c^*\on{Imm}(S^1,\mathbb R^2)\cong \mathcal D(S^1)^2$ 
is given by the integral operator $L\i$ which is
convolution by $F_n$ with respect to arc length $s$:
$$ \check{G}^{\on{imm},n}_c(h,k) = \iint_{S^1\times S^1} F_n(s_1-s_2).
\langle h(s_1),k(s_2)\rangle.ds_1.ds_2.$$

\subsection{Geodesics in the $G^{\on{imm},n}$-metric}\label{geodesicsGimm}
Differentiating the operator $D_s = \frac{1}{|c_\th|}\partial_\th$ with
respect to $c$ in the direction $m$ we get $-\frac{\langle m_\th, c_\th
\rangle}{|c_\th|^3} \partial _\th$, or $ -\langle D_s m, v \rangle D_s$.
Thus differentiating the big operator $L_{n,c}$ with respect to $c$ in the
direction $m$, we get:
\begin{align}
D_{(c,m)}L_{n,c}(h) &= (-1)^{n+1} A.\sum_{j=0}^{2n-1}  D_s^{j}
\langle  D_s(m), v\rangle D_s^{2n-j}(h)
\tag{3}\end{align}
Thus we have
\begin{align*}
&D_{(c,m)}G^{\on{imm},n}_c(h,k) =
\\&
= A.\int_{S^1} (-1)^{n+1} \sum_{j=0}^{2n-1}
\Big\langle  D_s^{j}
\langle  D_sm, v\rangle D_s^{2n-j}(h),k\Big\rangle ds
+\int_{S^1}\langle L_n(h),k\rangle \langle  D_sm, v\rangle  ds
\\&
= A.\int_{S^1} \sum_{j=1}^{2n-1}(-1)^{n+j+1}
\Big\langle \langle  D_sm, v\rangle  D_s^{2n-j}(h),
 D_s^{j}k\Big\rangle ds
+\int_{S^1}\langle h,k\rangle \langle  D_sm, v\rangle  ds
\\&
= \int_{S^1}
\bigg\langle m, A.\sum_{j=1}^{2n-1}(-1)^{n+j}
 D_s\bigg(\langle  D_s^{2n-j}h, D_s^{j}k\rangle v\bigg)
- D_s(\langle h,k \rangle v)
\bigg\rangle
 ds
\end{align*}
According to \ref{setting} we should rewrite this as
\begin{align*}
D_{(c,m)}G^{\on{imm},n}_c(h,k) &=
G^{\on{imm},n}_c(K^n_c(m,h), k)
=G^{\on{imm},n}_c\big(m, H^n_c(h,k)\big),
\end{align*}
and thus we find the two versions $K^n$ and $H^n$ of the $G^n$-gradient of
$c\mapsto G^{\on{imm},n}_c(h,k)$ are given by:
\begin{align*}
K^n_c(m,h) &= L_n\i\bigg((-1)^{n+1}A.\sum_{j=1}^{2n-1}
 D_s^{j} \langle  D_sm, v\rangle
   D_s^{2n-j}(h) +\langle D_sm, v\rangle h \bigg)
\tag{4}
\end{align*}
and by
\begin{align*}
&H^n_c(h,k) = L_n\i\bigg(A.\sum_{j=1}^{2n-1}(-1)^{n+j}
 D_s\Big(\langle  D_s^{2n-j}h,
 D_s^{j}k\rangle
 v\Big)- D_s(\langle h,k \rangle v)
\bigg)
\\
= & L_n\i\bigg(A.\sum_{j=1}^{2n-1}(-1)^{n+j} \langle  D_s^{2n-j+1}h,
D_s^{j}k\rangle v +A.\sum_{j=2}^{2i}(-1)^{n+j-1} \langle D_s^{2n-j+1}h,
D_s^{j}k\rangle v
\\& +A.\sum_{j=1}^{2n-1}(-1)^{n+j} \langle  D_s^{2n-j}h, D_s^{j}k\rangle
\ka_c n -\langle D_sh,k\rangle v -\langle h, D_sk\rangle v -\langle
h,k\rangle \ka_c n
\bigg) \\
=& L_n\i\bigg(
-\langle L_n(h),D_sk\rangle v
-\langle D_sh, L_n(k)\rangle v
-\langle h,k\rangle\ka(c) n \\
& \qquad  +A.\sum_{j=1}^{2n-1}(-1)^{n+j}
\langle  D_s^{2n-j}h,
 D_s^{j}k\rangle
\ka(c) n
\bigg)
\tag{5}\end{align*}
since
$ D_s\big( v\big)
  =\ka(c) n$.
By \ref{geodesic} the geodesic equation for the metric $G^n$ is
\begin{align*}
c_{tt}&= \tfrac12 H^n_c(c_t,c_t) - K^n_c(c_t,c_t).
\end{align*}
We expand it to get:
\begin{equation*}
\boxed{\quad
\begin{aligned}
L_n(c_{tt}) &= - \langle L_n(c_t),  D_s(c_t)\rangle v
-\frac{|c_t|^2\ka(c)}{2} n
-\langle D_s(c_t), v\rangle c_t
\\&\quad +\frac{A}{2}. \sum_{j=1}^{2n-1}(-1)^{n+j} \langle  D_s^{2n-j}c_t,
D_s^{j}c_t\rangle \ka(c) n
\\&\quad + (-1)^n A.\sum_{j=1}^{2n-1} D_s^{j}\left( \langle  D_s(c_t),
v\rangle D_s^{2n-j}(c_t)\right)
\end{aligned}
\quad}\tag{6}\end{equation*}
 From \thetag{3} we see that
\begin{align*}
(L_n(c_t))_t - L_n(c_{tt}) &=  dL_n(c)(c_t)(c_t) = (-1)^{n+1}
A.\sum_{j=0}^{2n-1} D_s^{j} \langle D_s(c_t),v\rangle D_s^{2n-j}(c_t).
\end{align*}
so that a more compact form of the geodesic equation of the metric $G^n$
is:
\begin{equation*}
\boxed{\quad
\begin{aligned}
(L_n(c_{t}))_t &=
- \langle L_n(c_t),  D_s(c_t)\rangle
 v
-\frac{|c_t|^2\ka(c)}{2} n
-\langle D_s(c_t), v\rangle L_n c_t
\\&\quad
+\frac{A}{2}.\sum_{j=1}^{2n-1}(-1)^{n+j}
\langle  D_s^{2n-j}c_t,
 D_s^{j}c_t\rangle
\ka(c) n
\end{aligned}
\quad}\tag{7}\end{equation*}
For $n=0$ this agrees with \cite{MM1},~4.1.2.

\subsection{Existence of geodesics}
\begin{thm*}\label{existenceGimm}
Let $n\ge 1$. For each $k\ge 2n+1$ the geodesic equation \ref{geodesicsGimm} \thetag{6}
has unique local solutions in the Sobolev space of
$H^{k}$-immersions. The solutions depend $C^\infty$ on $t$ and on the initial
conditions $c(0,\;.\;)$ and $c_t(0,\;.\;)$.  
The domain of existence (in $t$) is uniform in $k$ and thus this
also holds in $\on{Imm}(S^1,\mathbb R^2)$.
\end{thm*}

\begin{demo}{Proof}
We consider the geodesic equation as the flow equation of a smooth
($C^\infty$) vector field on the $H^2$-open set $U^k\x H^k(S^1,\mathbb R^2)$ 
in the Sobolev space $H^k(S^1,\mathbb R^2)\x H^k(S^1,\mathbb R^2)$ where
$U^k=\{c\in H^k:|c_\th|>0\}\subset H^k$ is $H^2$-open. 
To see that this works we will use the following facts: 
By the Sobolev inequality we have a
bounded linear embedding $H^{k}(S^1,\mathbb R^2)\subset C^m(S^1,\mathbb R^2)$ if
$k>m+\frac12$. The Sobolev space $H^k(S^1,\mathbb R)$ is a Banach algebra
under pointwise multiplication if $k>\frac12$. For any fixed smooth mapping
$f$ the mapping $u\mapsto f\o u$ is smooth $H^k\to H^k$ if $k>0$.
The mapping $(c,u)\mapsto L_{n,c}u$ is smooth $U\x H^k \to H^{k-2n}$ and is
a bibounded linear isomorphism $H^k\to H^{k-2n}$ for fixed $c$. This can be
seen as follows (see \ref{index-lemma} below): 
It is true if $c$ is parametrized by arclength (look at it
in the space of Fourier coefficients). The index is invariant under
continuous deformations of elliptic operators of fixed degree, 
so the index of $L_{n,c}$ is zero in general. But
$L_{n,c}$ is self-adjoint positive, so it is injective with vanishing
index, thus surjective. By the open mapping theorem it is then bibounded.
Moreover $(c,w)\mapsto L_{n,c}\i(w)$ is smooth $U^k\x H^{k-2n}\to H^k$ (by the
inverse function theorem on Banach spaces).
The mapping $(c,f)\mapsto D_s f = \frac1{|c_\th|}\p_\th f$
is smooth $H^k\x H^m\supset U\x H^m\to H^{m-1}$ for $k\ge m$, and is linear in $f$.
Let us write 
$D_c f=D_s f$ just for the remainder of this proof to stress the dependence
on $c$.
We have $v=D_c c$ and $n=J D_c c$. The mapping $c\mapsto \ka(c)$ is smooth
on the $H^2$-open set $\{c:|c_\th|>0\}\subset H^k$ into $H^{k-2}$.
Keeping all this in mind we now write the geodesic equation as follows: 
\begin{align*}
c_t &= u =:X_1(c,u)
\\
u_t &= L_{n,c}\i\Big(- \langle L_{n,c}(u),  D_c(u)\rangle D_c(c)
  -\frac{|c_t|^2\ka(c)}{2} JD_c(c) -\langle D_c(u), D_c c\rangle u
\\&\qquad\qquad 
+\frac{A}{2}. \sum_{j=1}^{2n-1}(-1)^{n+j} \langle  D_c^{2n-j}u, 
  D_c^{j}u\rangle \ka(c) JD_c(c)
\\&\qquad\qquad 
+ (-1)^n A.\sum_{j=1}^{2n-1} D_c^{j}\left( \langle  D_c(u), D_c(c)
  \rangle D_c^{2n-j}(u)\right)\Big)
\\&
=: X_2(c,u)
\end{align*}
Now a term by term investigation of this shows that 
the expression in the brackets is smooth $U^k\x H^k\to H^{k-2n}$ since
$k-2n\ge 1>\frac12$. The operator $L_{n,c}\i$ then takes it smoothly back
to $H^k$. So the vector field $X=(X_1,X_2)$ is smooth on $U^k\x H^k$.
Thus the flow $\on{Fl}^k$ exists on $H^k$ and is smooth in $t$ and the initial 
conditions for fixed $k$.

Now we consider smooth initial conditions $c_0=c(0,\;.\;)$ and
$u_0=c_t(0,\;.\;)=u(0,\;.\;)$ in $C^\infty(S^1,\mathbb R^2)$. Suppose the 
trajectory $\on{Fl}^k_t(c_0,u_0)$ of
$X$ through these intial conditions in $H^k$ maximally exists for $t\in
(-a_k,b_k)$, and the trajectory $\on{Fl}^{k+1}_t(c_0,u_0)$ in $H^{k+1}$ 
maximally exists for $t\in(-a_{k+1},b_{k+1})$ with $b_{k+1}<b_k$. By
uniqueness we have $\on{Fl}^{k+1}_t(c_0,u_0)=\on{Fl}^{k}_t(c_0,u_0)$ for
$t\in (-a_{k+1,}b_{k+1})$. We now apply $\p_\th$ to the equation
$u_t=X_2(c,u)=L_{n,c}\i(\,\dots\,)$, 
note that the commutator $[\p_\th,L_{n,c}\i]$ is a pseudo differential
operator of order $-2n$ again, and write $w=\p_\th u$. We obtain 
$w_t=\p_\th u_t= L_{n,c}\i\p_\th(\,\dots\,) + [\p_\th,L_{n,c}\i](\,\dots\,)$.
In the term $\p_\th(\,\dots\,)$ we consider now only the terms $\p_\th^{2n+1}u$ and
rename them $\p_\th^{2n}w$.
Then we get an 
equation $w_t(t,\th)=\tilde X_2(t,w(t,\th))$ which is inhomogeneous bounded linear
in $w\in H^k$ with
coefficients bounded linear operators on $H^k$ 
which are $C^\infty$ functions of $c, u \in H^k$. 
These we already know on the intervall $(-a_k,b_k)$. 
This equation 
therefore has a solution $w(t,\;.\;)$ for all $t$ for which the
coefficients exists, thus for all $t\in (a_k,b_k)$. The limit 
$\lim_{t\nearrow b_{k+1}} w(t,\;.\;)$ exists in $H^k$ 
and by continuity it equals $D_c(u)$ in $H^k$
at $t=b_{k+1}$. Thus the $H^{k+1}$-flow was not
maximal and can be continued. So $(-a_{k+1},b_{k+1})=(a_k,b_k)$.  
We can iterate this and conclude that the flow of $X$ exists in
$\bigcap_{m\ge k} H^{m}= C^\infty$. 
\qed\end{demo}

\subsection{The conserved momenta of $G^{\on{imm},n}$}
According to \ref{momentum} the following
momenta are preserved along any geodesic
$t\mapsto c(t,\;.\;)$:
$$
\boxed{
\begin{aligned}
&\langle c_\th,L_{n,c}(c_t)\rangle |c_\th(\th)|\;\in\X(S^1)
&\quad&\text{reparametrization momentum}
\\
&\int_{S^1}L_{n,c}(c_t)\,ds=\int_{S^1}c_t\,ds\;\in \mathbb R^2
&\quad&\text{linear momentum}
\\
&\int_{S^1}\langle Jc,L_{n,c}(c_t)\rangle \,ds\in \mathbb R
  &\quad&\text{angular momentum}
\end{aligned}
}
$$

\subsection{Horizontality for $G^{\on{imm},n}$}
\label{horizontal}
$h\in T_c\on{Imm}(S^1,\mathbb R^2)$ is
$G^{\on{imm},n}_c$-orthogonal to the
$\on{Diff}(S^1)$-orbit through $c$ if and only if
$$
0 =G^{\on{imm},n}_c(h,\ze_X(c)) = G^{\on{imm},n}_c(h,c_\th.X)
=\int_{S^1}X.\langle L_{n,c}(h),c_\th \rangle \,ds
$$
for all $X\in\X(S^1)$. So the $G^{\on{imm},n}$-normal bundle is given by
\begin{equation*}
\mathcal N^n_c = \{h\in C^\infty(S,\mathbb R^2):
  \langle L_{n,c}(h), v \rangle = 0\}.
\end{equation*}
The $G^{\on{imm},n}$-orthonormal projection $T_c\on{Imm}\to \mathcal N^n_c$, denoted
by $h\mapsto h^\bot=h^{\bot,G^n}$ and the complementary projection
$h\mapsto h^\top\in T_c(c\o\on{Diff}(S^1))$ are determined as follows:
$$
h^\top = X(h).v \quad\text{  where }
\langle L_{n,c}(h),v \rangle = \langle L_{n,c}(X(h).v),v \rangle
$$
Thus we are led to consider the linear differential operators associated to
$L_{n.c}$
\begin{align*}
L^\top_c, L^\bot_c&: C^\infty(S^1) \to C^\infty(S^1),
\\
L^\top_c(f)&=\langle L_{n,c}(f.v),v \rangle
  =\langle L_{n,c}(f.n),n \rangle,
\\
L^\bot_c(f)&=\langle L_{n,c}(f.v),n \rangle
  = -\langle L_{n,c}(f.n),v \rangle.
\end{align*}
The operator $L^\top_c$  is of order $2n$ and also unbounded, self-adjoint
and positive on $L^2(S^1,|c_\th|\,d\th)$ since
\begin{align*}
\int_{S^1}L^\top_c(f)g ds
&=\int_{S^1}\langle L_{n,c}(fv), v\rangle g ds
\\&
=\int_{S^1}\langle fv, L_{n,c}(gv)\rangle  ds
=\int_{S^1}f L^\top_c(g) ds,
\\
\int_{S^1}L^\top_c(f)f ds
&=\int_{S^1}\langle fv, L_{n,c}(fv)\rangle  ds \quad
 > 0 \text{  if }f\ne0.
\end{align*}
In particular, $L^\top_c$ is injective.
$L^\bot_c$, on the other hand is of order $2n-1$ and a similar argument
shows it is skew-adjoint.
For example, if $n=1$, then one finds that:
\begin{align*}
L^\top_c &= - A.D_s^2 +(1+A.\ka^2).I \\
L^\bot_c &= -2A.\ka. D_s - A.D_s(\ka).I
\end{align*}

\begin{lem*}\label{index-lemma}
The operator $L^\top_c:C^\infty(S^1)\to C^\infty(S^1)$ is invertible.
\end{lem*}

\begin{demo}{Proof}
This is because its index vanishes, by the
following argument: The index is invariant under continuous deformations of
elliptic operators of degree $2n$. The operator
$$
L^\top_c(f) =(-1)^n\frac{A}{|c_\th|^{2n}}\p_\th^{2n}(f)
+ \text{  lower order terms}
$$
is homotopic to $(1+(-1)^n\p_\th^{2n})(f)$ and thus has the same index
which
is zero since the operator $1+(-1)^n\p_\th^{2n}$ is invertible. This can be
seen by expanding in Fourier series where the latter operator is given by
$(\hat f(m))\mapsto ((1+m^{2n})\hat f(m))$, a linear isomorphism of the
space of rapidly decreasing sequences.
Since $B^\top_c$ is injective, it is also surjective.
\qed\end{demo}

To go back and forth between the `natural' horizontal space of vector
fields $a.n$ and the $G^{\on{imm},n}$-horizontal vector fields $\{h \mid
\langle Lh,v \rangle =0\}$, we only need to use these operators and the
inverse of $L^\top$. Thus, given $a$, we want to find $b$ and $f$ such that
$L(an+bv)=fn$, so that $an+bv$ is $H^n$-horizontal. But this implies that
$$ L^\bot(a) = -\langle L(an),v \rangle = \langle L(bv),v \rangle =
L^\top(b).$$
Thus if we define the operator $C_c:C^\infty(S^1)\to C^\infty(S^1)$ by
$$ C_c :=  (L^\top_c)\i \circ L^\bot_c, $$
we get a pseudo-differential operator of order -1 (which is an integral
operator), so that $a.n + C(a).v$ is always $H^{\text{imm},n}$-horizontal. In
particular, the restriction of the metric $G^{\on{imm},n}$ to horizontal
vector fields $h_i = a_i.n + b_i.v$ can be computed like this:
\begin{align*}
G^{\on{imm},n}_c(h_1,h_2) &= \int_{S^1} \langle Lh_1, h_2 \rangle .ds \\
&= \int_{S^1} \langle L(a_1.n + b_1.v), n \rangle. a_2.ds \\
&= \int_{S^1} \left( L^\top(a_1) + L^\bot(b_1) \right).a_2.ds \\
&= \int_{S^1} \left( L^\top + L^\bot \circ C \right) a_1.a_2.ds.
\end{align*}
Thus the metric restricted to horizontal vector fields is given by the
pseudo differential operator $L^{\on{red}} = L^\top + L^\bot \circ
(L^\top)\i \circ L^\bot.$ On the quotient space $B_i$, 
if we identify its tangent space at $C$ with the space 
of normal vector fields $a.n$, then:
$$\boxed{G^{imm,n}_C(a_1,a_2) = \int_C (L^\top + L^\bot \circ
(L^\top)\i \circ L^\bot)a_1 \cdot a_2 \cdot ds}$$
Now, although this operator may be hard to analyze, its inverse, the metric
on the cotangent space to $B_i$, is simple. The tangent space to $B_i$ at a
curve $C$ is canonically the quotient of that of $\on{Imm}(S^1, \mathbb R^2)$ 
at a parametrization $c$ of $C$, modulo the subspace of multiples of
$v$. Hence the cotangent space to $B_i$ at $C$ injects into that of
$\on{Imm}(S^1, \mathbb R^2)$ at $c$ with image the linear functionals that
vanish on $v$. In terms of the dual basis $\check{v}, \check{n}$, these are
multiples of $\check{n}$. On the smooth cotangent space 
$C^\infty(S^1,\mathbb R^2)\cong G^0_c(T_c\on{Imm}(S^1,\mathbb R^2))
\subset T_c^*\on{Imm}(S^1,\mathbb R^2)\cong \mathcal D(S^1)^2$ 
the dual metric is given by convolution with the elementary kernel
$K_n$ which is a simple sum of exponentials. Thus we need only restrict
this kernel to multiples $a(s).\check{n}_c(s)$ to obtain the dual metric on
$B_i$. The result is that:
$$ \check{G}^n_c(a_1,a_2) = \iint_{S^1\times S^1} K_n(s_1-s_2). \langle
n_c(s_1), n_c(s_2) \rangle .a_1(s_1). a_2(s_2). ds_1 ds_2.$$

\subsection{Horizontal geodesics}
The normal bundle $\mathcal N_c$ mentioned in \ref{horizontal} is well
defined
and is a smooth vector subbundle of the tangent bundle.
But $\on{Imm}(S^1,\mathbb R^2)\to B_i(S^1,\mathbb
R^2)=\on{Imm}/\on{Diff}(S^1)$
is {\it not} a principal bundle and thus there are no principal connections, but we
can prove the main consequence, the existence of horizontal paths,
directly:

\begin{prop*}
For any smooth path $c$ in $\on{Imm}(S^1,\mathbb R^2)$ there exists a
smooth
path $\ph$ in $\on{Diff}(S^1)$ with $\ph(0,\;.\;)=\on{Id}_{S^1}$ depending
smoothly on $c$ such that
the path $e$ given by $e(t,\th)=c(t,\ph(t,\th))$ is horizontal:
$\langle L_{n,e}(e_t),e_\th\rangle=0$.
\end{prop*}

\begin{demo}{Proof}
Writing $D_c$ instead of $D_s$ we note that 
$D_{c\o\ph}(f\o \ph)=\frac{(f_\th\o\ph)\ph_\th}{|c_\th\o\ph|.|\ph_\th|}
=(D_c(f))\o\ph$ for $\ph\in\on{Diff}^+(S^1)$. So we have 
$L_{n,c\o\ph}(f\o\ph)=(L_{n,c}f)\o\ph$.

Let us write $e=c\o \ph$ for $e(t,\th)=c(t,\ph(t,\th))$, etc.
We look for $\ph$
as the integral curve of a time dependent vector field $\xi(t,\th)$ on
$S^1$, given by $\ph_t=\xi\o \ph$.
We want the following expression to vanish:
\begin{align*}
\langle L_{n,c\o\ph}(\p_t(c\o\ph)),\p_\th(c\o\ph) \rangle
&=\langle L_{n,c\o\ph}(c_t\o\ph + (c_\th\o\ph)\,\ph_t),(c_\th\o\ph)\,\ph_\th
\rangle
\\&
=\langle L_{n,c}(c_t)\o\ph + L_{n,c}(c_\th.\xi)\o\ph,c_\th\o\ph
\rangle\ph_\th
\\&
=\bigl((\langle L_{n,c}(c_t),c_\th\rangle +\langle L_{n,c}(\xi.c_\th),
c_\th\rangle)\o\ph\bigr)\,\ph_\th.
\end{align*}
Using the time dependent vector field
$\xi=-\frac1{|c_\th|}(L^\top_c)\i(\langle L_{n,c}(c_t),v\rangle)$
and its flow $\ph$ achieves this.
\qed\end{demo}

If we write
$$c_t =na + vb=\Big(n,v\Big)\binom{a}{b}$$
then we can expand the condition for horizontality as follows:
\begin{align*}
D_s(c_t)
&=\big(D_sa+\ka(c) b\big)n
+\big(D_sb-\ka(c) a\big)v.
\\&
=(n,v)
  \begin{pmatrix} D_s & \ka \\
                  -\ka   &  D_s\end{pmatrix}
  \binom{a}{b}
\\
L_n^c(c_t)
&= c_t +(-1)^n A (n,v)
   \begin{pmatrix} D_s & \ka \\
                  -\ka   &  D_s\end{pmatrix}^{2n}
   \binom{a}{b}
\\&
=c_t + (-1)^n A (n,v)
\begin{pmatrix} D_s^2-\ka^2 & D_s\ka +\ka D_s\\
  -D_s\ka -\ka D_s  & D_s^2-\ka^2  \end{pmatrix}^n
   \binom{a}{b}
\end{align*}
so that horizontality becomes
\begin{align*}
0&= \langle L_{n,c}(c_t),v\rangle
=\langle c_t,v\rangle + (-1)^n (0,1)
\begin{pmatrix} D_s^2-\ka^2 & D_s\ka +\ka D_s\\
  -D_s\ka -\ka D_s & D_s^2-\ka^2 \end{pmatrix}^{n}
   \binom{a}{b}
\end{align*}

We may specialize the general geodesic equation to horizontal paths and
then take the $v$ and $n$ parts of the geodesic equation. For a horizontal
path we may write $L_{n,c}(c_t)=\tilde an$ for
$\tilde a(t,\th)=\langle L_{n,c}(c_t),n \rangle$.
The $v$ part of the equation turns out to vanish identically and then $n$
part gives us (because $n_t$ is a multiple of $v$):
$$ \boxed{\quad
\tilde a_t
= -\frac{|c_t|^2\ka(c)}{2}
-\langle D_s c_t,v\rangle
\tilde a
+\frac{\ka(c)}2 \sum_{j=1}^{2n-1}(-1)^{n+j}
\langle D_s^{2n-j}c_t,D_s^{j}c_t\rangle \quad}
$$
Note that applying \ref{existenceGimm} with horizontal initial vectors
gives us local existence and uniqueness for solutions of this horizontal
geodesic equation.


\subsection{A Lipschitz bound for arclength in $G^{\on{imm},n}$}
We apply the inequality of Cauchy-Schwarz to the
derivative
of the length function
$\ell(c)=\int|c_\th|d\th$ along a path
$t\mapsto c(t,\;.\;)$:
\begin{align*}
\p_t\ell(c)&=d\ell(c)(c_t)
=\int_{S^1} \frac{\langle c_{t\th}, c_\th\rangle}{|c_\th|}\, d\th
=\int_{S^1} \langle  D_s(c_t), v\rangle\, ds
\\&
\le \Big(\int_{S^1}| D_s(c_t)|^2 ds\Big)^\frac12
    \cdot \Big(\int_{S^1}1^2 ds\Big)^\frac12
\le \sqrt{\ell(c)}\frac1A \|c_t\|_{G^1},
\\&
\le \sqrt{\ell(c)}C(A,n)\|c_t\|_{G^n},
\\
\p_t\sqrt{\ell(c)} &= \frac{\p_t\ell(c)}{2\sqrt{\ell(c)}}
  \le\frac{C(A,n)}2 \|c_t\|_{G^n}.
\end{align*}
Thus we get
\begin{align*}
|\sqrt{\ell(c(1,\;.\;))}-\sqrt{\ell(c(0,\;.\;))}|
&\le \int_0^1|\p_t\sqrt{\ell(c)}|\,dt
\le \frac{C(A,n)}2\int_0^1\|c_t\|_{G^n}\,dt 
\\&
= \frac{C(A,n)}2 L_{G^n}(c).
\end{align*}
Taking the infinimum of this over all paths $t\mapsto c(t,\;.\;)$ from
$c_0$ to $c_1$ we see that for $n\ge 1$ we have the Lipschitz estimate:
$$
|\sqrt{\ell(c_1)}-\sqrt{\ell(c_0)}|\le
\frac12\on{dist}_{G^n}^{\on{Imm}}(c_1,c_0) 
$$
Since we have $L_{G^n}^{\text{hor}}(c)\le L_{G^n}(c)$ with equality for
horizontal curves we also have:
$$
\boxed{\text{{\it If} }n\ge 1, \quad |\sqrt{\ell(C_1)}-\sqrt{\ell(C_0)}|\le
\frac12\on{dist}_{G^n}^{B_i}(C_1,C_0) }
$$

\subsection{Scale invariant immersion Sobolev metrics}
Let us mention in passing that we may use the length of the curve to modify
the immersion Sobolev metric so that it becomes scale invariant:
\begin{align*}
G^{\text{imm,scal},n}_c(h,k) &= \int_{S^1}\big(\ell(c)^{-3}\langle h,k \rangle 
+ \ell(c)^{2n-3} A \langle D_s^n(h),D_s^n(k) \rangle\big)\,ds
\\&
= \int_{S^1}\big\langle(\ell(c)^{-3}+(-1)^n \ell(c)^{2n-3}A D_s^{2n})h,k\big \rangle\,ds
\end{align*}
This metric can easily be analyzed using the methods described above. In
particular we note that the geodesic equation on $\on{Imm}(S^1,\mathbb
R^2)$ for this metric is built in a similar way than that for
$G^{\text{imm},n}$ and that the existence theorem in \ref{existenceGimm}
holds for it. Note the conserved momenta along a geodesic $t\mapsto
c(t,\;.\;)$ are:
\begin{align*}
& \frac{1}{\ell(c)^{3}}\int_{S^1}c_t\,ds 
  + (-1)^n \ell(c)^{2n-3}A\int_{S^1} D_s^{2n}(c_t)\,ds 
\\&\quad 
= \frac{1}{\ell(c)^{3}}\int_{S^1}c_t\,ds 
\in \mathbb R^2
&\quad &\text{linear momentum}
\\&
\frac{1}{\ell(c)^{3}}\int_{S^1}\langle Jc,c_t\rangle\,ds 
+(-1)^n \ell(c)^{2n-3}A\int_{S^1}\langle Jc, D_s^{2n}(c_t)\rangle\,ds
&\quad &\text{angular momentum}
\\&
\frac{1}{\ell(c)^{3}}\int_{S^1}\langle c,c_t\rangle\,ds
+(-1)^n \ell(c)^{2n-3}A\int_{S^1}\langle c, D_s^{2n}(c_t)\rangle\,ds
&\quad &\text{scaling momentum}
\end{align*}
As in the work of Trouv\'e and Younes \cite{YT, Y}, we may consider the
following variant.
\begin{align*}
G^{\text{imm,scal},n,\infty}_c(h,k) 
&=\lim_{A\to\infty}\frac1A 
\int_{S^1}\big\langle(\ell(c)^{-3}+(-1)^n \ell(c)^{2n-3}A D_s^{2n})h,k\big \rangle\,ds
\\&
=(-1)^n \ell(c)^{2n-3}\int_{S^1}\big\langle D_s^{2n}h,k\big \rangle\,ds
\end{align*}
It is degenerate with kernel the constant tangent vectors. The interesting
fact is that the scaling momentum for $G^{\text{imm,scal},1,\infty}$ is given by 
$$
-\frac1{\ell(c)}\int_{S^1}\langle c,D_s^2(c_t) \rangle ds = \p_t \log\ell(c).
$$

\section{Sobolev metrics on $\on{Diff}(\mathbb R^2)$ and on its
quotients}\label{G^{diff,n}}

\subsection{The metric on $\on{Diff}(\mathbb R^2)$.}
\label{diffR2}
We consider the regular Lie group $\on{Diff}(\mathbb R^2)$ which is either
the group $\on{Diff}_c(\mathbb R^2)$ of all diffeomorphisms with compact
supports of $\mathbb R^2$ or the group $\on{Diff}_{\mathcal S}(\mathbb
R^2)$ of all diffeomorphisms which decrease rapidly to the identity. The
Lie algebra is $\X(\mathbb R^2)$, by which we denote either the Lie algebra
$\X_c(\mathbb R^2)$ of vector fields with compact support or the Lie
algebra $\X_{\mathcal S}(\mathbb R^2)$ of rapidly decreasing vector fields,
with the negative of the usual Lie bracket. For any $n \ge 0$, we equip
$\on{Diff}(\mathbb R^2)$ with the right invariant weak Riemannian metric
$G^{\on{Diff},n}$given by the Sobolev $H^n$-inner product on $\X_c(\mathbb
R^2)$.
\begin{align*}
H^n(X,Y)&=\underset{\substack{0\le i,j\le n\\i+j\le n}}{\sum}
\frac{A^{i+j}n!}{i!j!(n-i-j)!}
 \int_{\mathbb R^2}  \langle \p_{x^1}^i\p_{x^2}^j X,\p_{x^1}^i\p_{x^2}^j Y
 \rangle\,dx
\\&
=\underset{\substack{0\le i,j\le n\\i+j\le n}}{\sum}
(-A)^{i+j} \frac{n!}{i!j!(n-i-j)!}
\int_{\mathbb R^2}
  \langle \p_{x^1}^{2i}\p_{x^2}^{2j} X,Y \rangle\,dx
\\&
=\int_{\mathbb R^2} \langle L X,Y \rangle\,dx \qquad\text{  where }
\\
L &= L_{A,n} = (1 - A\De)^n, \qquad \De = \p_{x^1}^2+\p_{x^2}^2.
\end{align*}
(We will write out the full subscript of $L$ only where it helps clarify
the meaning.) The completion of $\X_c(\mathbb R^2)$ is the Sobolev space
$H^n(\mathbb R^2)^2$. With the usual $L^2$-inner product we can identify
the dual of $H^n(\mathbb R^2)^2$ with $H^{-n}(\mathbb R^2)^2$ (in the space
of tempered distributions). Note that the operator $L:H^n(\mathbb R^2)^2\to
H^{-n}(\mathbb R^2)^2$ is a bounded linear operator. On $L^2(\mathbb R^2)$
the operator $L$ is unbounded selfadjoint and positive. In terms of Fourier
transform we have $\widehat{L_{A,n}u}(\xi)=(1+A|\xi|^2)^n\hat u$. Let
$F_{A,n}$ in the space of tempered distributions $\mathcal S'(\mathbb R^2)$ 
be the fundamental solution (or
Green's function: note that we use the letter `F' for `fundamental' because
`G' has been used as the metric) of $L_{A,n}$ satisfying
$L_{A,n}(F_{A,n})=\de_0$ which is
given by
$$ F_{A,n}(x) =\frac1{2\pi}\int_{\mathbb R^2}e^{i\langle x,\xi\rangle}
\frac{1}{(1+A|\xi|^2)^n}\,d\xi. $$
The functions $F_{A,n}$ are given by the classical modified Bessel
functions $K_r$ (in the notation, e.g., of Abramowitz and Stegun 
\cite{AbramowitzStegun} or of
Matlab) by the
formula:
$$  F_{A,n}(x) =
\frac{1}{2^n\pi (n-1)! A}\, .\left(\frac{|x|}{\sqrt{A}}\right)^{n-1} 
K_{n-1}\left(\frac{|x|}{\sqrt{A}}\right).$$
and it satisfies $(L\i u)(x)=\int_{\mathbb R^2}F(x-y)u(y)\,dy$ for each
tempered distribution $u$. The function $F_{A,n}$ is $C^{n-1}$ except that
$F_{A,1}$ has a log-pole at zero. At infinity, $F_{A,n}(x)$ is
asymptotically a constant times $x^{n-3/2}e^{-x}$: these facts plus much
much more can be found in \cite{AbramowitzStegun}.

\subsection{Strong conservation of momentum and `EPDiff'}
What is the form of the conservation of momentum for a geodesic $\ph(t)$ in
this metric, that is to say, a flow $x \mapsto \ph(x,t)$ on $\mathbb R^2$?
We need to work out $\on{Ad}_\ph^*$ first. Using the definition, we see:
\begin{align*}
\int_{\mathbb R^2} \langle LX, \on{Ad}_\ph^*(Y) \rangle
&
:=\int_{\mathbb R^2} \langle L\on{Ad}_\ph(X), Y
\rangle = \int_{\mathbb R^2} \langle (d\ph . X)\circ \ph\i,LY \rangle \\
&= \int_{\mathbb R^2} \det(d\ph) \langle d\ph . X, LY \circ \ph \rangle
= \int_{\mathbb R^2} \langle X, \det(d\ph). d\ph^T . (LY \circ \ph)
\rangle
\end{align*}
hence:
$$ \on{Ad}_\ph^*(Y) = L\i \left(\det(d\ph) . d\ph^T . (LY \circ \ph)
\right).$$
Now the conservation of momentum for geodesics $\ph(t)$ 
of right invariant metrics on
groups says that:
$$ L\i \left( \det(d\ph)(t) . d\ph(t)^t . \left( L(\frac{\p \ph}{\p t}
\circ \ph\i ) \circ \ph \right)\right)$$
is independent of $t$. This can be put in a much more transparent form.
First, $L$ doesn't depend on $t$, so we cross out the outer $L\i$. Now let
$v(t) = \frac{\p \ph}{\p t} \circ \ph\i \in \X(\mathbb R^2)$ be the tangent
vector to the geodesic. Let $u(t) = Lv(t)$, so that:
$$ \det(d\ph)(t) . d\ph(t)^t . (u(t) \circ \ph(t))$$
is independent of $t$.
We should {\it not} think of $u(t)$ as a vector field on $\mathbb R^2$:
this is because we want $\langle u,v \rangle$ to make invariant sense in
any coordinates whatsoever. This means we should think of $u$ as expanding
to the differential form:
$$ \om(t) = (u_1.dx^1 + u_2.dx^2)\otimes \mu$$
where $\mu = dx^1\wedge dx^2$, the area form. But then:
$$ \ph(t)^*(\om(t)) = \langle d\ph^t . (u \circ \ph(t)), dx
\rangle \otimes \det(d\ph)(t).\mu$$
so conservation of momentum says simply:
$$ \boxed{\quad \ph(t)^*\om(t) \text{ is independent of } t \quad}$$
This motivates calling $\om(t)$ the momentum of the geodesic flow.
As we mentioned above, conservation of momentum for a Riemannian metric on
a group is very strong and is an integrated form of the geodesic equation.
To see this, we need only take the differential form of this conservation
law. $v(t)$ is the infinitesimal flow, so the infinitesimal form of the
conservation is:
$$ \frac{\p}{\p t} \om(t) + \mathcal L_{v(t)}(\om(t)) = 0$$
where $\mathcal L_{v(t)}$ is the Lie derivative. We can expand this term
by term:
\begin{align*}
\mathcal L_{v(t)}(u_i) &= \sum_j v^j.\frac{\p u_i}{\p x^j} \\
\mathcal L_{v(t)}(dx^i) &= dv^i = \sum_j \frac{\p v^i}{\p x^j}.dx^j\\
\mathcal L_{v(t)}(\mu) &= \on{div}v(t)\mu \\
\mathcal L_{v(t)}(\om(t)) &=  \left(\sum_{i,j}\left( v^j.\frac{\p u_i}{\p
x^j}.dx^i + u^j.\frac{\p v^j}{\p x^i}.dx^i\right) + \on{div}v.\sum_i u_i
dx^i\right) \otimes \mu.
\end{align*}
The resulting differential equation for geodesics has been named {\it
EPDiff}:
$$
\boxed{\quad\begin{aligned}
 &v=\frac{\p \ph}{\p t} \o \ph\i,\qquad u = L(v)
 \\
 &\frac{\p u_i}{\p t} 
 + \sum_j\left( v^j.\frac{\p u_i}{\p x^j} 
 + u^j.\frac{\p v^j}{\p x^i}\right) + \on{div}v. u_i = 0.
 \end{aligned}\quad}
$$
Note that this is a special case of the general equation of Arnold: $\p_t
u=-\on{ad}(u)^*u$ for geodesics on any Lie group in any right (or left)
invariant metric. The name `EPDiff' was coined by Holm and Marsden and
stands for `Euler-Poincar\'e', although it takes a leap of faith to see it
in the reference they give to Poincar\'e.

\subsection{The quotient metric on $\on{Emb}(S^1,\mathbb R^2)$}
\label{Gdiffn}
We now consider the quotient mapping $\on{Diff}(\mathbb R^2) \to
\on{Emb}(S^1,\mathbb R^2)$ given by $\ph \mapsto \ph \circ i$ as in the
section \ref{nmb:1}. Since this identifies $\on{Emb}(S^1, \mathbb R^2)$
with a right coset space
of $\on{Diff}(\mathbb R^2)$, and since the metric $G_{\on{diff}}^n$ is
right invariant, we can put a quotient metric on $\on{Emb}(S^1, \mathbb
R^2)$ for which this map is a Riemannian submersion. Our next step is to
identify this metric. Let $\ph \in \on{Diff}(\mathbb R^2)$ and let $c = \ph
\circ i \in \on{Emb}(S^1, \mathbb R^2)$. The fibre of this map through
$\ph$ is the coset
$$
\ph.\on{Diff}^0(S^1,\mathbb R^2) = \{\ps \mid  \ps \circ c \equiv c \}.\ph.
$$
whose tangent space is (the right translate by $\ph$ of) 
the vector space of vector fields $X\in \X(\mathbb R^2)$ 
with $X \circ c \equiv 0$. This is the vertical subspace. Thus the
horizontal subspace is
$$\left\{Y \left| \int_{\mathbb R^2} \langle LY,X \rangle dx = 0,
  \text{ if } X\circ c \equiv 0 \right. \right\}. $$
If we want $Y\in \X(\mathbb R^2)$ then the horizontal subspace is 0. But we
can also search for $Y$ in a bigger space of vector fields on $\mathbb
R^2$. What we need is that $LY=c_*(p(\th).ds)$, where $p$ is a function
from $S^1$ to $\mathbb R^2$ and $ds$ is arc-length measure supported on
$C$. To make $c_*(p(\th).ds)$ pair with smooth vector fields $\X(\mathbb
R^2)$ in a coordinate invariant way, we should interpret the values of $p$
as 1-forms. Solving for $Y$, we have:
\begin{align*}
Y(x) &= \int_{S^1}F(x-c(\th)).p(\th)ds
\end{align*}
(where, to make $Y$ a vector field, the values of $p$ are now interpreted
as vectors, using the standard metric on $\mathbb R^2$ to convert 1-forms
to vectors). Because $F$ is not $C^\infty$, we have a case here where the
horizontal subspace is not given by $C^\infty$ vector fields. However, we
can still identify the set of vector fields in this horizontal subspace
which map bijectively to the $C^\infty$ tangent space to $\on{Emb}(S^1,
\mathbb R^2)$ at $c$. Mapped to $T_c\on{Emb}(S^1,\mathbb R^2)$, the above
$Y$ goes to:
\begin{align*}
(Y\o c)(\th) &= \int_{S^1}F(c(\th)-c(\th_1)).p(\th_1).|c'(\th_1)|d\th_1
\\&
=: (F_c \ast p)(\th) \quad\text{  where } \tag{1}
\\
F_c(\th_1,\th_2) &= F(c(\th_1)-c(\th_2))
= \frac1{2\pi}\int_{\mathbb R^2}e^{i\langle c(\th_1)-c(\th_2),\xi\rangle}
\frac{1}{(1+A|\xi|^2)^n}\,d\xi.
\end{align*}
Note that here, convolution on $S^1$ uses the metric
$L^2(S^1,|c'(\th|d\th)$ and it defines a self-adjoint operator for this
Hilbert space. Moreover, it is covariant with respect to change in
parametrization:
$$ F_{c\circ\ph} \ast (f \circ \ph) = (F_c \ast f) \circ \ph.$$

What are the properties of the kernel $F_c$? From the properties of $F$, we
see that $F_c$ is $C^{n-1}$ kernel (except for log poles at the diagonal
when $n=1$). It is also a pseudo-differential operator of order $-2n+1$ on
$S^1$. To see that let us assume for the moment that each function of $\th$
is a periodic function on $\mathbb R$. Then
\begin{align*}
c(\th_1)-c(\th_2)&=\int_0^1 c_\th(\th_2+t(\th_1-\th_2))dt.(\th_1-\th_2)
=:\tilde c(\th_1,\th_2)(\th_1-\th_2)
\\
F_c(\th_1,\th_2) &
= \frac1{2\pi}\int_{\mathbb R^2}e^{i(\th_1-\th_2)\langle \tilde c(\th_1,\th_2),\xi\rangle}
\frac{1}{(1+A|\xi|^2)^n}\,d\xi
\\&
= \frac1{2\pi}\int_{\mathbb R}e^{i(\th_1-\th_2)\et_1}
\Big(\int_{\mathbb R} \frac{|\tilde c(\th_1,\th_2)|^{-2}}
{(1+\frac{A}{|\tilde c(\th_1,\th_2)|^{2}}(|\et_1|^2+|\et_2|^2))^n}\,d\et_2\Big)\,d\et_1
\\&
=: \frac1{2\pi}\int_{\mathbb R}e^{i(\th_1-\th_2)\et_1}
\tilde F_c(\th_1,\th_2,\et_1)\,d\et_1
\end{align*}
where we changed variables as $\et_1 = \langle \tilde c(\th_1,\th_2),\xi\rangle$
and $\et_2 = \langle J \tilde c(\th_1,\th_2),\xi \rangle$. 
So we see that $F_c(\th_1,\th_2)$ is an elliptic pseudo differential
operator kernel of degree $-2n+1$ (the loss comes from integrating with
respect to $\et_2$). The symbol $\tilde F_c$ is real and positive, so the
operator $p\mapsto F_c * p$ is self-adjoint and positive. Thus it is
injective, and by an index argument similar to the one in \ref{index-lemma}
it is invertible. The inverse
operator to the integral operator $F_c$ is a pseudo-differential operator
$L_c$ of order $2n-1$ given by the distribution kernel $L_c(\th,\th_1)$
which satisfies
\begin{align*}
L_c \ast F_c \ast f &= F_c \ast L_c \ast f = f 
\\\tag{2} 
L_{c\circ \ph} \ast (h \circ \ph) 
&= ((L_c \ast h)\circ \ph) \quad \text{for all } \ph \in \on{Diff}^+(S^1)
\end{align*}

If we write $h = Y \circ c$, then we want to express the horizontal lift
$Y$ in terms of $h$ and write $Y_h$ for it. The set of all these $Y_h$
spans the horizontal subspace which maps isomorphically to
$T_c\on{Emb}(S^1,\mathbb R^2)$. Now:
$$ h = Y \circ c = \left(F \ast (c_* (p.ds))\right)\circ c = F_c \ast p.$$
Therefore, using the inverse operator, we get $p = L_c \ast h$ and:
\begin{align*} 
Y_h &= F\ast (c_* (p.ds)) = F\ast(c_* ((L_c\ast h).ds)) \quad \text{or}\\
Y_h(x) &= \int_{S^1}F(x-c(\th))\int_{S^1}L_c(\th,\th_1)h(\th_1)
  |c'(\th_1)|d\th_1|c'(\th)|d\th
\end{align*}
and $LY_h = c_*((L_c \ast h).ds)$. 
Thus we can finally write down the quotient metric
\begin{align*}
G^{\on{diff},n}_c&(h,k) = \int_{\mathbb R^2} \langle LY_h,Y_k \rangle dx \\
& = \int_{S^1} \Big\langle
L_c \ast h(\th),\int_{S^1}F(c(\th)-c(\th_1))
\int_{S^1}L_c(\th_1,\th_2)k(\th_2)ds_2  \,ds_1 \Big\rangle ds
\tag{3}\\
&= \int_{S^1} \langle L_c \ast h(\th), k(\th) \rangle ds 
= \iint_{S^1\times S^1} L_c(\th,\th_1) \langle h(\th_1), k(\th) \rangle ds_1\, ds.
\end{align*}
The dual metric on the associated smooth cotangent space $L_c* C^\infty(S^1,\mathbb R^2)$ is similarly:
$$ \check G^{\on{diff},n}_c(p,q) 
= \iint_{S^1\times S^1} F_c(\th,\th_1) \langle p(\th_1), q(\th) \rangle ds_1\, ds.$$


\subsection{The geodesic equation on $\on{Emb}(S^1, \mathbb R^2)$ via
conservation of momentum}\label{geodA}
A quite convincing but not rigorous derivation of this equation can be
given using the fact that under a submersion, geodesics on the quotient
space are the projections of those geodesics on the total space which are
horizontal at one and hence every point. In our case, the geodesics on
$\on{Diff}(\mathbb R^2)$ can be characterized by the strong conservation of
momentum we found above: $\ph(t)^*\om(t)$ is independent of $t$. If $X(t)$
is the tangent vector to the geodesic, i.e.\ the velocity $X(t) =
\partial_t \ph \circ \ph\i (t)$, then $\om(t)$ is just
$LX(t)=c_*(p(\th,t).ds) = c_*(p(\th,t).|c_\th(\th,t)|.d\th)$ considered as
a measure valued 1-form instead of a vector field.

When we pass to the quotient $\on{Emb}(S^1,\mathbb R^2)$, a horizontal
geodesic of diffeomorphisms $\ph(t)$ with $\ph(0) = \text{identity}$ gives
a geodesic path of embeddings $c(\th,t) = \ph(t) \circ c(0,\th)$. For these
geodesic equations, it will be most convenient to take as the momentum the
1-form ${\tilde p}(\th,t)=p(\th,t).|c_\th(\th,t)|$, the measure factor $d\th$ being
constant along the flow. We must take the velocity to be the horizontal
vector field $X(t) = F \ast c(\cdot,t)_*({\tilde p}(\th,t).d\th)$. For this to be
the velocity of the path of maps $c$, we must have $c_t(\th,t)=X(c(\th),t)$
because the global vector field $X$ must extend $c_t$.
To pair ${\tilde p}$ and $c_t$, we regard ${\tilde p}$ as a 1-form along $c$
(the area factor having been replaced by the measure $d\th$ supported on
$C$).
The geodesic equation must be the differential form of the conservation
equation:
$$\boxed{ \ph(t)^* {\tilde p}(\cdot,t) \text{ is independent of } t.}$$
More explicitly, if $d_x$ stands for differentiating with respect to the
spatial coordinates $x,y$, then this means:
$$ d_x\ph(t)^T |_{c(\th,t)} {\tilde p}(\th,t) = \text{ cnst.}$$
We differentiate this with respect to $t$, using the identity:
$$ \partial_t d_x\ph(t) = d_x(\ph_t(t)) = d_x(X \circ \ph(t)) 
= (d_x(X)\circ \ph(t)) \cdot d_x\ph(t),$$
we get
$$ 0 = d_x \ph(t)^T \cdot \left((d_x(X)^T \circ c(\th,t)) \cdot
{\tilde p}(\th,t) + {\tilde p}_t(\th,t) \right).$$
Writing this out and putting the discussion together, we get the following
form for the geodesic equation on $\on{Emb}(S^1, \mathbb R^2)$:
$$\begin{aligned}
c_t(\th,t) &= X(t) \circ c(\th,t)\\
{\tilde p}_t(\th,t) &= -\on{grad} X^t (c(\th,t),t)\cdot {\tilde p}(\th,t)\\
X(t) &= F \ast c(\cdot,t)_* \left({\tilde p}(\th,t).d\th\right)
\end{aligned}$$
Note that $X$ is a vector field on the plane: these are not closed
equations if we restrict $X$ to the curves. The gradient of $X$ requires
that we know the normal derivative of $X$ to the curves. Alternatively, we
may introduce a second {\it vector-valued} kernel on $S^1$ depending on $c$
by:
$$F'_c(\th_1,\th_2) = \on{grad} F(c(\th_1)-c(\th_2)).$$
Then the geodesic equations may be written:
$$\boxed{\begin{aligned}
c_t(\th,t) &= (F_c \ast {\tilde p})(\th,t) \\
{\tilde p}_t(\th,t) &= -\langle {\tilde p}(\th,t),(F'_c \ast {\tilde p})(\th,t) \rangle.
\end{aligned}}$$
where, in the second formula, the dot product is between the two $\tilde
p$'s and the vector value is given by $F'_c$.

The problem with this approach is that we need to enlarge the space
$\on{Diff}(\mathbb R^2)$ to include diffeomorphisms which are not
$C^\infty$ along some $C^\infty$ curve but have a mild singularity normal
to the curve. Then we would have to develop differential geometry and the
theory of geodesics on this space, etc. It seems more straightforward to
outline the direct derivation of the above geodesic equation, along the
lines used above.

\subsection{The geodesic equation on $\on{Emb}(S^1, \mathbb R^2)$, direct
approach}\label{geodesicsGdiff}
The space of invertible pseudo differential operators on a compact
manifold is a regular Lie group (see 
\cite{ARS2}
), so we can use
the usual formula $d(A\i)=-A\i.dA.A\i$ for computing the derivative of
$L_c$ with respect to $c$. Note that we have a simple expression for
$D_{c,h}F_c$, namely
\begin{align*}
D_{c,h}F_c(\th_1,\th_2) &= dF(c(\th_1)-c(\th_2))(h(\th_1)-h(\th_2))
=\langle F'_c(\th_1,\th_2),h(\th_1)-h(\th_2) \rangle \\
\text{hence}\qquad \qquad &\\
D_{c,\ell}L_c(\th_1,\th_2) &= -\int_{(S^1)^2}
L_c(\th_1,\th_3)\,D_{c,h}F_c(\th_3,\th_4)\,L_c(\th_4,\th_2)\,d\th_3\,d\th_4\\
&= -\int_{(S^1)^2} L_c(\th_1,\th_3)\, \langle
(F'_c(\th_3,\th_4),\ell(\th_3)\rangle \,L_c(\th_4,\th_2)\,d\th_3\,d\th_4 \\
&\quad +\int_{(S^1)^2} L_c(\th_1,\th_4)\, \langle
(F'_c(\th_4,\th_3),\ell(\th_3)\rangle \,L_c(\th_3,\th_2)\,d\th_3\,d\th_4
\end{align*}
We can now differentiate the metric where $\th=(\th_1,\th_2,\dots,\th_n)$ is the
variable on $(S^1)^n$:
\begin{align*}
&D_{c,\ell}G^{\text{diff},n}(h,k)
= \int_{(S^1)^2} D_{c,\ell} L_c(\th_1,\th_2) \langle h(\th_2),
k(\th_1)\rangle d\th
\\&
= \int_{(S^1)^4} \Big\langle 
  -L_c(\th_1,\th_3) F'_c(\th_3,\th_4) L_c(\th_4,\th_2)
\\&\qquad\qquad
+ L_c(\th_1,\th_4) F'_c(\th_4,\th_3)\,L_c(\th_3,\th_2),
  \ell(\th_3)\Big\rangle\langle h(\th_2), k(\th_1)\rangle\,d\th 
\end{align*}
We have to write this in the form 
$$
D_{c,\ell}G^{\text{diff},n}_c(h,k) = G^{\text{diff},n}_c(\ell,H_c(h,k)) 
= G^{\text{diff},n}_c(K_c(\ell,h),k)
$$
For $H_c$ we use $\de(\th_5-\th_3)=\int
L_c(\th_5,\th_6)F_c(\th_6,\th_3)d\th_6 = (L_c * F_c)(\th_5,\th_3)$ as follows: 
\begin{align*}
&D_{c,\ell}G^{\text{diff},n}(h,k)  
=\int_{(S^1)^6} L_c(\th_5,\th_6)\Big\langle \ell(\th_5),
  \Big(-L_c(\th_1,\th_3) F'_c(\th_3,\th_4) L_c(\th_4,\th_2)
\\&\qquad\qquad
+ L_c(\th_1,\th_4) F'_c(\th_4,\th_3)\,L_c(\th_3,\th_2)\Big) F_c(\th_6,\th_3)
  \langle h(\th_2), k(\th_1)\rangle\Big\rangle\,d\th 
\end{align*}
Thus 
\begin{align*}
H_c(h,k)(\th_0) &= \int_{(S^1)^4}
\Big(-L_c(\th_1,\th_3) F'_c(\th_3,\th_4) L_c(\th_4,\th_2)
\\&\qquad\qquad
+ L_c(\th_1,\th_4) F'_c(\th_4,\th_3)\,L_c(\th_3,\th_2)\Big) F_c(\th_0,\th_3)
  \langle h(\th_2), k(\th_1)\rangle\,d\th
\end{align*}
Similarly we get
\begin{align*}
&D_{c,\ell}G^{\text{diff},n}(h,k)  
= \int_{(S^1)^6} L_c(\th_6,\th_5) \Big\langle  F_c(\th_1,\th_6)
  \big\langle -L_c(\th_1,\th_3)  F'_c(\th_3,\th_4) L_c(\th_4,\th_2)
\\&\qquad\qquad
+ L_c(\th_1,\th_4) F'_c(\th_4,\th_3)\,L_c(\th_3,\th_2),
  \ell(\th_3)\big\rangle h(\th_2), k(\th_5)\Big\rangle\,d\th 
\end{align*}
so that 
\begin{align*}
&K_c(\ell,h)(\th_0) 
=\\&
= \int_{(S^1)^2}  
  \Big(-\langle F'_c(\th_0,\th_1) L_c(\th_1,\th_2),\ell(\th_0)\rangle
+ \langle F'_c(\th_0,\th_1)\,L_c(\th_1,\th_2),
  \ell(\th_1)\rangle\Big) h(\th_2)\,d\th
\end{align*}
By \ref{geodesic} the geodesic equation is given by 
\begin{align*}
c_{tt}(\th_0) &= \tfrac12 H_c(c_t,c_t)(\th_0) - K_c(c_t,c_t)(\th_0)
\end{align*}
Let us rewrite the geodesic equation in terms of $L_c*c_t$. 
We have (suppressing the variable $t$ and collecting all terms)
\begin{align*}
(L_c*c_t)_t(\th_0) &= \int_{S^1} D_{c,c_t}L_c(\th_0,\th_1)c_t(\th_1)\,d\th_1
  + L_c*c_{tt}
\\&
= \tfrac12 \int_{S^1}
\Big(F'_c(\th_1,\th_0) - F'_c(\th_0,\th_1)\Big)
  \langle L_c*c_t(\th_0), L_c*c_t(\th_1)\rangle\,d\th_1
\end{align*}
Since the kernel $F$ is an even function we get the same geodesic equation as above for the momentum ${\tilde p}(\th,t)=L_c*c_t=p(\th,t).|c_\th|$:
\begin{equation*}
\begin{aligned}
{\tilde p}_t(\th_0) &= - \int_{S^1} F_c'(\th_0,\th_1)\langle {\tilde p}(\th_0), {\tilde p}(\th_1) \rangle\,d\th_1
\end{aligned}\tag{1}
\end{equation*}

\subsection{Existence of geodesics}
\begin{thm*}\label{existenceG^diff}
Let $n\ge 1$. For each $k> 2n-\frac12$ the geodesic equation
\ref{geodesicsGdiff} \thetag{1}
has unique local solutions in the Sobolev space of
$H^{k}$-embeddings. The solutions are $C^\infty$ in $t$ and in the initial
conditions $c(0,\;.\;)$ and $c_t(0,\;.\;)$.  
The domain of existence (in $t$) is uniform in $k$ and thus this
also holds in $\on{Emb}(S^1,\mathbb R^2)$. 
\end{thm*}

An even stronger theorem, proving {\it global} existence on the level of
$H^k$-diffeo\-mor\-phisms on $\mathbb R^2$, has been proved by 
\cite{T, YT1, YT2}.

\begin{demo}{Proof} Let $c\in H^{k}$.
We begin by checking that $F'_c$ is a pseudo differential operator kernel
of order $-2n+2$ as we did for $F_c$ in \ref{Gdiffn}.
\begin{align*}
c(\th_1)-c(\th_2)&
=:\tilde c(\th_1,\th_2)(\th_1-\th_2)
\\
\on{grad}F(x)  
&= \frac1{2\pi}\int_{\mathbb R^2}e^{i\langle x,\xi\rangle}
\frac{J\xi}{(1+A|\xi|^2)^n}\,d\xi
\\
F'_c(\th_1,\th_2) &
= \frac1{2\pi}\int_{\mathbb R^2}e^{i(\th_1-\th_2)\langle \tilde c(\th_1,\th_2),\xi\rangle}
\frac{J\xi}{(1+A|\xi|^2)^n}\,d\xi
\\&
= \frac1{2\pi}\int_{\mathbb R}e^{i(\th_1-\th_2)\et_1}
\Big(\int_{\mathbb R} \frac{|\tilde c(\th_1,\th_2)|^{-3}.J\et}
{(1+\frac{A}{|\tilde c(\th_1,\th_2)|^{2}}(|\et_1|^2+|\et_2|^2))^n}\,d\et_2\Big)\,d\et_1
\\&
=: \frac1{2\pi}\int_{\mathbb R}e^{i(\th_1-\th_2)\et_1}
\tilde F_c(\th_1,\th_2,\et_1)\,d\et_1
\end{align*}
where we changed variables as $\et_1 = \langle \tilde c(\th_1,\th_2),\xi\rangle$
and $\et_2 = \langle J \tilde c(\th_1,\th_2),\xi \rangle$. 
So we see that $F'_c(\th_1,\th_2)$ is an elliptic pseudo differential
operator kernel of degree $-2n+2$ 
(the loss of $1$ comes from integrating with
respect to $\et_2$).  
We write the geodesic equation in the following way:
\begin{align*}
c_t &= F_c*q =: Y_1(c,q)
\\
q_t&= \big\langle q, F'_c* q \big\rangle 
= \int F'_c(\;.\;,\th)\langle q(\th),q(\;.\;) \rangle \,d\th =: Y_2(c,q)
\end{align*}
We start with $c\in H^k$ where $k>2n-\frac12$, in the $H^2$-open set
$U^k:=\{c:|c_\th|>0\}\subset H^k$. Then $q=L_c*c_t\in H^{k-2n+1}$ and 
$F'_c*q\in H^{k-1}\subset H^{k-2n+1}$. 
By the Banach algebra property of the Sobolev space $H^{k-2n+1}$ 
the expression (with missuse of notation)
$Y_2(c,q)=\langle q,F'_c*q \rangle\in H^{k-2n+1}$.
Since the kernel $F$ is not smooth only at 0, all appearing pseudo
differential operators kernels are $C^\infty$ off the diagonal, thus are
smooth mappings in $c$ with values in the space of operators between the
relevant Sobolev spaces. Let us make this more precise. We claim that 
$c\mapsto F'_c*(\;.\;)\in L(H^k, H^{k+2n-2})$ is $C^\infty$. Since the
Sobolev spaces are convenient, we can (a) use the smooth uniform
boundedness theorem \cite{KM},~5.18, so that it suffices to check that for
each fixed $q\in H^k$ the mapping $c\mapsto F'_c*q$ is smooth into
$H^{k+2n-2}$. Moreover, by
\cite{KM},~2.14 it suffices (b) to check that this is weakly smooth: Using
the $L^2$-duality between $H^{k+2n-2}$ and $H^{-k-2n+2}$ it suffices to
check, that for each $p\in H^{-k-2n+2}$ the expression
\begin{align*}
&\int p(\th_1)(F'_c* q)(\th_1)\,d\th_1 =
\\&
=\iint \frac{p(\th_1)}{2\pi}\int_{\mathbb R}e^{i(\th_1-\th_2)\et_1}
\Big(\int_{\mathbb R} \frac{|\tilde c(\th_1,\th_2)|^{-3}.J\et}
{(1+\frac{A}{|\tilde
c(\th_1,\th_2)|^{2}}(|\et_1|^2+|\et_2|^2))^n}\,d\et_2\Big)\,d\et_1
q(\th_2)\,d\th_1\,d\th_2
\end{align*}
is a smooth mapping $\on{Emb}(S^1,\mathbb R^2)\to \mathbb R^2$. For that we
may assume that $c$ depends on a further smooth variable $s$. Convergence
of this integral depends on the highest order term in the asymptotic
expansion in $\et$, which does not change if we differentiate with respect
to $s$.

Thus the geodesic equation is the flow equation of
a smooth vector field $Y=(Y_1,Y_2)$ on $U^k\x H^{k-2n+1}$. We thus have
local existence and uniqueness of the flow $\on{Fl}^k$ on $U^k\x H^{k-2n+1}$.

Now we consider smooth initial conditions $c_0=c(0,\;.\;)$ and
$q_0=q(0,\;.\;)=(L_c*c_t)(0,\;.\;)$ in $C^\infty(S^1,\mathbb R^2)$. 
Suppose the trajectory $\on{Fl}^k_t(c_0,q_0)$ of
$Y$ through these intial conditions in $U^k\x H^{k+1-2n}$ 
maximally exists for $t\in (-a_k,b_k)$, and the trajectory 
$\on{Fl}^{k+1}_t(c_0,u_0)$ in $U^{k+1}\x H^{k+2-2n}$ 
maximally exists for $t\in(-a_{k+1},b_{k+1})$ with $b_{k+1}<b_k$. By
uniqueness we have $\on{Fl}^{k+1}_t(c_0,u_0)=\on{Fl}^{k}_t(c_0,u_0)$ for
$t\in (-a_{k+1,}b_{k+1})$. We now apply $\p_\th$ to the equation
$q_t=Y_2(c,q)$, note that the commutator
$q\mapsto [F'_c,\p_\th]*q=\p_th(F'_c*q)-F'_c*(\p_\th q)$ 
is a pseudo differential operator
of order $-2n+2$ again, and obtain 
\begin{align*}
\p_\th q_t &= \int [F'_c,\p_\th](\;.\;,\th)\langle q(\th),q(\;.\;) \rangle \,d\th
+ \int F'_c(\;.\;,\th)\langle \p_\th q(\th),q(\;.\;) \rangle \,d\th
\\&\quad
+ \int F'_c(\;.\;,\th)\langle  q(\th), \p_\th q(\;.\;) \rangle \,d\th
\end{align*}
which is an inhomogeneous linear equation for $w=\p_\th q$ in $U^k\x
H^{k+1-2n}$.
By the variation of constant method one sees that the solution 
$w$ exists in $H^k$ for as long as $(c,q)$ exists in 
$U^k\x H^{k+1-2n}$, i.e., for all $t\in (-a_k,b_k)$. By continuity we can
conclude that $w=\p_\th q$ is the derivative in $H^{k+2-2n}$ for
$t=b_{k+1}$, and thus the domain of definition was not maximal. Iterating
this argument we can conclude that the solution $(c,q)$ lies in $C^\infty$
for $t\in (-a_k,b_k)$.
\qed\end{demo}

\subsection{Horizontality for $G^{\on{diff},n}$}
The tangent  vector $h\in T_c\on{Emb}(S^1,\mathbb R^2)$ is
$G^{\on{diff},n}_c$-orthogonal to the
$\on{Diff}(S^1)$-orbit through $c$ if and only if
$$ 0 =G^{\on{diff},n}_c(h,\ze_X(c))=
\int_{(S^1)^2}L_c(\th_1,\th_2)\langle h(\th_2),c_\th(\th_1) \rangle X(\th_1) \,ds_1\,ds_2 $$
for all $X\in\X(S^1)$. So the $G^{\on{diff},n}$-normal bundle is given by
\begin{equation*}
\mathcal N^{\on{diff},n}_c = \{h\in C^\infty(S^1,\mathbb R^2):
  \langle L_c * h, v \rangle = 0\}.
\end{equation*}
Working exactly as in section 4, we want to split any tangent vector into
vertical and horizontal parts as $h = h^\top + h^\bot$ where $h^\top =
X(h).v$ for $X(h)\in \X(S^1)$ and where $h^\bot$ is horizontal, $\langle
L_c * h^\bot, v \rangle = 0$. Then $\langle L_c * h, v \rangle = \langle
L_c * (X(h)v),v \rangle$ and we are led to consider the following
operators:
\begin{align*}
&L_c^\top, L_c^\bot: C^\infty(S^1) \to  C^\infty(S^1),\\
& L_c^\top (f) = \langle L_c*(f.v),v \rangle = \langle L_c*(f.n),n \rangle, \\
& L_c^\bot (f) = \langle L_c*(f.v),n \rangle = -\langle L_c*(f.n),v \rangle.
\end{align*}
The pseudo differential operator $L_c^\top$ is unbounded, selfadjoint and
positive on $L^2(S^1,d\th)$ since we have 
\begin{align*}
\int_{S^1}L_c^\top(f).f\,d\th &= 
\int_{(S^1)^2}\langle L_c(\th_1,\th_2)f(\th_2)v(\th_2),f(\th_1).v(\th_1)\rangle\,d\th 
= \|f.v\|^2_{G^{\on{diff},n}} > 0.
\end{align*}
Thus $L_c^\top$ is injective and by an index argument as in
\ref{index-lemma} the operator $L_c^\top$ is invertible. Moreover, the
operator $L_c^\bot$ is skew-adjoint. To go back and forth between the
natural horizontal space of vector fields $a.n$ and the
$G^{\on{diff},n}$-horizontal vectors, we have to find $b$ such that $L_c *
(a.n+b.v) = f.n$ for some $f$. But then
$$
L_c^\bot(a)= -\langle L_c*(a.n), v \rangle = \langle L_c*(b.v),v \rangle 
= L_c^\top(b)\quad\text{  thus }\quad b = (L_c^\top)\i L_c^\bot(a).
$$
Thus $a.n + (L_c^\top)\i L_c^\bot(a).v$ is always 
$G^{\on{diff},n}$-horizontal and is the 
horizontal projection of $a.n+b.v$ for any $b$.

\begin{prop*}
For any smooth path $c$ in $\on{Imm}(S^1,\mathbb R^2)$ there exists a
smooth path $\ph$ in $\on{Diff}(S^1)$ with $\ph(0,\;.\;)=\on{Id}_{S^1}$
depending
smoothly on $c$ such that the path $e$ given by $e(t,\th)=c(t,\ph(t,\th))$
is $G^{\on{diff},n}$-horizontal: $\langle L_c * e_t,e_\th\rangle=0$.
\end{prop*}

\begin{demo}{Proof}
Let us write $e=c\o \ph$ for $e(t,\th)=c(t,\ph(t,\th))$, etc.
We look for $\ph$
as the integral curve of a time dependent vector field $\xi(t,\th)$ on
$S^1$, given by $\ph_t=\xi\o \ph$.
We want the following expression to vanish. In its computation the
equivariance of $L_c$ under $\ph\in\on{Diff}^+(S^1)$ from
\ref{Gdiffn}\thetag2 will
play an important role.
\begin{align*}
&\big\langle L_{c\o \ph}*(\p_t(c\o\ph)),\p_\th(c\o\ph) \big\rangle
=\langle  L_{c\o\ph}*(c_t\o\ph + (c_\th\o\ph)\,\ph_t),(c_\th\o\ph)\,\ph_\th
\rangle
\\&
=\big\langle((L_c * c_t)\o\ph)\ + ((L_c*(c_\th.\xi))\o\ph),
(c_\th\o\ph)\ph_\th\big\rangle
\\&
=\big((\langle L_c*c_t,c_\th\rangle 
  +\langle L_c*(\xi.c_\th),c_\th\rangle)\o\ph\bigr)\,\ph_\th.
\end{align*}
Using the time dependent vector field
$\xi=-(L^\top_c)\i\langle L_c*c_t,c_\th\rangle$
and its flow $\ph$ achieves this.
\qed\end{demo}

To write the quotient metric on $B_e$, we want to lift normal vector fields
$a.n$ to a curve $C$ to horizontal vector fields on $\on{Emb}(S^1,\mathbb
R^2)$. Substituting $h = a.n+(L_c^\top)\i L_c^\bot(a).v,\,\, k =
b.n+(L^\top)\i L^\bot(b).v$ in \ref{Gdiffn}\thetag3, we get as above:
$$ G^{\on{diff},n}_C (a,b) = \int_C \left(L_c^\top+ L_c^\bot (L_c^\top)\i
L_c^\bot\right)(a).b ds.$$
The dual metric on the cotangent space is just the restriction of the dual
metric on $\on{Emb}(S^1,\mathbb R^2)$ to the cotangent space to $B_e$ and
is much simpler. We simply set $p=f.n, \,\, q=g.n$ and get:
$$ \check G^{\on{diff},n}_C(f,g) = \iint_{C^2} F(x(s)-x(s_1)). \langle
n(s), n(s_1) \rangle .f(s)g(t).dsds_1$$
where $x(s) \in \mathbb R^2$ stands for the point in the plane with arc
length coordinate $s$ and $F$ is the Bessel kernel. Since these are dual
inner products, we find that the two operators, (a) convolution with the
kernel $F(x(s)-x(s_1)).\langle n(s), n(s_1) \rangle$ and (b) $L^\top_c +
L^\bot_c (L^\top_c)^{-1} L^\bot_c$ are inverses of each other.

\subsection{The geodesic equation on $B_e$ via conservation of momentum}
The simplest way to find the geodesic equation on $B_e$ is again to
specialize the general rule $\ph(t)^*\om(t) = \text{ cnst.}$ to the
horizontal geodesics.
Now horizontal in the present context, that is for $B_e$, requires more of
the momentum $\om(t)$. As well as being given by $c_*(p(s).ds)$, we require
the 1-form $p$ to kill the tangent vectors $v$ to the curve. If we identify
1-forms and vectors using the Euclidean metric, then we may say simply
$p(s)=a(s).n$, where $a$ is a scalar function on $C$. But note that if you
take the momentum as $c_*(a(s)n(s)ds)$ and integrate it against a vector
field $X$, then you find:
$$ \langle X, c_*(a(s)n(s)ds) \rangle = \int_C a(s) \langle X, n(s) \rangle
ds = \int_C a(s).i_X(dx \wedge dy)$$
where $i_X$ is the `interior product' or contraction with $X$ taking a
2-form to a 1-form. Noting that 1-forms can be integrated along curves
without using any metric, we see that the 2-form along $c$ defined by
$\{a(s).(dx\wedge dy)_{c(s)}\}$ can be naturally paired with vector fields
so it defines a canonical measure valued 1-form. Therefore, the momentum
for horizontal geodesics can be identified with this 2-form.

If $\ph(x,t)$ is a horizontal geodesic in $\on{Diff}(\mathbb R^2)$, then
the curves $C_t = \on{image}(c(\cdot,t))$ are given by $C_t=\ph(C_0,t)$ and
the momentum is given by $a(\th,t).(dx \wedge dy)$, where $c(\th,t)$
parametrizes the curves $C_t$. Note that in order to differentiate $a$ with
respect to $t$, we need to assign parameters on the curves $C_t$
simultaneously. We do this in the same way we did for almost local metrics:
assume $c_\th$ is a multiple of the normal vector $n_C$. But $\th_0 \mapsto
\ph(\th_0, t)$ gives a second map from $C_0$ to $C_t$: in terms of $\th$,
assume this is $\th = \bar \ph(\th_0,t)$. Then the conservation of momentum
means simply:
$$ a(\bar \ph(\th_0,t),t).\on{det}(D_x\ph)(c(\th_0,0),t) \text{ is independent of } t.$$
Let $X$ be the global vector field giving this geodesic, so that $\ph_t = X
\circ \ph$. Note that $\bar\ph_t = (\langle X\circ c,v \rangle /|c_\th|)
\circ \bar \ph.$ Using this fact, we can differentiate the displayed
identity. Recalling the definition of the flow from its momentum and the
identifying $T_C B_e$ with normal vector fields along $C$, we get the full
equations for the geodesic:
\begin{align*}
C_t &= \langle X,n \rangle \cdot n\\
a_t &= -\langle X, v \rangle D_s(a) - \on{div}(X).a \\
X &= F \ast c_*(a(s)n(s)ds)
\end{align*}
Note, as in the geodesic equations in \ref{geodA}, that we must use $F$ to
extend $X$ to the whole plane. In this case, we only need (a) the normal
component of $X$ along $C$, (b) its tangential component along $C$ and (c)
the divergence of $X$ along $C$. These are obtained by convolving $a(s)$
with the kernels (which we give now in terms of arc-length):
\begin{align*}
F_c^{nn}(s_1,s_2) &= F(c(s_1)-c(s_2)) \langle n(s_1), n(s_2) \rangle ds_2\\
F_c^{vn}(s_1,s_2) &= F(c(s_1)-c(s_2)) \langle v(s_1), n(s_2) \rangle ds_2\\
F_c^{\on{div}}(s_1,s_2) &= \langle \on{grad}F(c(s_1)-c(s_2)), n(s_1) \rangle ds_2
\end{align*}
Then the geodesic equations become:
\begin{align*}
C_t &= (F_c^{nn} \ast a).n\\
a_t &= -(F_c^{vn} \ast a) D_s(a) - (F_c^{\on{div}} \ast a).a 
\end{align*}

Alternately, we may specialize the geodesic equation in \ref{geodA} to
horizontal paths. Then the $v$ part vanishes identically and the $n$ part
gives the last equation above. We omit this calculation.

%
%

\section{Examples}\label{examples}

\subsection{The Geodesic of concentric circles}
\label{concentric1}
All the metrics that we have studied are invariant under the motion group,
thus
the 1-dimensional submanifold of $B_e$ consisting of all concentric circles
centered at the origin is the fixed point set of the group of rotations
around the center. Therefore it is a geodesic in all our metrics. It is
given by the map $c(t,\th)= r(t)e^{i\th}$. Then $c_\th=ire^{i\th}$,
$v_c=ie^{i\th}$, $n_c=-e^{i\th}$, $\ell(c)=2\pi r(t)$, $\ka_c=\frac1{r}$
and $c_t=r_te^{i\th}=-r_t.n_c$.

The parametrization $r(t)$ can be determined by requiring constant speed
$\si$, i.e.\ if the metric is $G(h,k)$, then we require $G_c(c_t,c_t) =
r_t^2 G_c(n_c,n_c) = \si^2$, which leads to 
$\sqrt{G_c(c_t,c_t)}\,dt= \pm \sqrt{G_c(n_c,n_c)}\,dr$. 
To determine when the geodesic is complete as
$ r \rightarrow 0$ and $ r\rightarrow \infty$, we merely need to look at
its length which is given by:
$$ \int_0^\infty \sqrt{G_c(c_t,c_t)} dt=\int_0^\infty \sqrt{G_c(n_c,n_c)} dr,$$
and we need to ask whether this integral converges or diverges at its two limits.
Let's consider this case by case.

\noindent
{\bf  The metric $G^\Ph$}:
The geodesic is determined by the equation:
$$ G^\Ph(c_t,c_t)=2\pi r \cdot\Ph\left(2\pi r(t), \frac{1}{r(t)}\right)\cdot r_t^2 
= \si^2.$$
Differentiating this with respect to $t$ leads to the geodesic equation in
the standard form $r_{tt}=r_t^2f(r)$. It is easily checked that all three
invariant momentum mappings vanish: the reparameterization, linear and
angular momentum.

\begin{thm*} If $\Phi(2\pi r, 1/r) \approx C_1 r^a \text{ (resp. } C_2
r^b)$ as $r \rightarrow 0 \text{ (resp. } \infty)$, then the geodesic of
concentric circles is complete for $r \rightarrow 0$ if and only $a \le -3$
and is complete for $r \rightarrow \infty$ if and only if $b \ge -3$. In
particular, for $\ph = \ell^k$, we find $k=a=b$ and the geodesic is given
by $r(t) = \text{cnst.}t^{2/(k+3}$. For the scale invariant case
$\Ph(\ell,\ka)=\frac{4\pi^2}{\ell^3}+\frac{\ka^2}{\ell}$, we find $a=b=-3$,
the geodesic is given by $r(t) = e^{\sqrt{2}\si t}$ and is complete.
Moreover, in this case, the scaling momentum $\frac{2r_t}{r}$ is constant
in $t$ along the geodesic.
\end{thm*}
The proof is straightforward.

\noindent{\bf The metric $G^{\text{imm},n}$}
Recall from \ref{H^n-geodesic} the operator 
$L_{n,c} = I + (-1)^nA.D_s^{2n}$. For $c(t,\th)=r(t)e^{i\th}$ 
\ref{concentric1} we have
$$
L_{n,c}(c_t)=\big(1+(-1)^n\frac{A}{r^{2n}}\p_\th^{2n}\big)(r_te^{i\th})
=r_t\big(1+\frac{A}{r^{2n}}\big)e^{i\th}
$$
which is still normal to $c_\th$. So $t\mapsto c(t,\;.\;)$ is a horizontal
path for any choice of $r(t)$. Thus its speed is the square root of:
$$ G^{\on{imm},n}(c_t,c_t) = 2\pi r \cdot \left(1+\frac{A}{r^{2n}}\right) \cdot r_t^2 =\si^2.$$
For $n=1$ this is the same as the identity for the metric with
$\Ph(\ell,\ka)=1+A\ka^2$ which was computed in \cite{MM1},~5.1. An
explanation of this phenomenon is in \cite{MM1},~3.2.

\begin{thm*} The geodesic of concentric cirles is complete in the
$G^{\on{imm},n}$ metric if $n \ge 2$. For $n=1$, it is incomplete as $r
\rightarrow 0$ but complete if $r \rightarrow \infty$.
\end{thm*}

\noindent{\bf The metric $G^{\on{diff},n}$}
To evaluate the norm of a path of concentric circles, we now need to find
the vector field $X$ on $\mathbb R^2$ gotten by convolving the Bessel
kernel with the unit normal vector field along a circle. Using circular
symmetry, we find that:
\begin{align*} &X(x,y) = f(r)\left(\frac{x}{r}, \frac{y}{r}\right)r_t\\
&\left(I-A(\partial_{rr}+\frac{1}{r} \partial_r)\right)^n f =0 
\text{ except on the circle } r=r_0 \\
& f \in C^{2n-2} \text{ everywhere}, f(r_0)=1
\end{align*}
For $n=1$, we can solve this and the result is the vector field on $\mathbb
R^2$given by the Bessel functions $I_1$ and $K_1$:
$$X(x,y) = \begin{cases} \frac{I_1(r/\sqrt{A})}{I_1(r_0/\sqrt{A})} \text{ if } r \le r_0 \\ 
\frac{K_1(r/\sqrt{A})}{K_1(r_0/\sqrt{A})} \text{ if } r \ge r_0 \end{cases}
$$
Using the fact that the Wronskian of $I_1,K_1$ is $1/r$, we find:
$$ G^{\on{diff},1}(r_t n, r_t n) = \int \langle (I-A\triangle) X,X\rangle r_t^2 
= \frac {2\pi r_t^2}{K_1(r/\sqrt{A}).I_1(r/\sqrt{A})}.$$
Using the asymptotic laws for Bessel functions, one finds that the geodesic
of concentric circles has finite length to $r=0$ but infinite length to
$r=\infty$.

For $n>1$, it gets harder to solve for $X$. But lower bounds are not hard:
\begin{align*}
G^{\on{diff},n}(n,n) 
&= \inf_{X,\langle X,n\rangle \equiv 1 \text{ on } C_r} 
\int \langle (I-A\triangle)^nX,X \rangle \\
&\ge A^n.\underset{\left(
\substack{
{X,\langle X,n\rangle \equiv 1 \text{ on } C_r} \\ 
{X \rightarrow 0, \text{ when } x \rightarrow \infty}}
\right)}
{\inf} \int \langle \triangle^n(X), X \rangle \underset{\text{def}}{=}
M(r)
\end{align*}
Then $M(r)$ scales with $r$: $M(r)=M(1)/r^{2n-2}$, hence the length of the
path when the radius shrinks to 0 is bounded below by $\int_0 dr/r^{n-1}$
which is infinite if $n>1$. On the other hand, the metric $G^{\on{diff},n}$
dominates the metric $G^{\on{diff},1}$ so the length of the path when the
radius grows to infinity is always infinite. Thus:
\begin{thm*} The geodesic of concentric cirles is complete in the
$G^{\on{diff},n}$ metric if $n \ge 2$. For $n=1$, it is incomplete as $r
\rightarrow 0$ but complete if $r \rightarrow \infty$.
\end{thm*}

\subsection{Unit balls in five metrics at a `cigar'-like shape}
It is useful to get a sense of how our various metrics differ. One way to
do this is to take one simple shape $C$ and examine the unit balls in the
tangent space $T_C B_e$ for various metrics. All of our metrics (except the
simple $L^2$ metric) involve a constant $A$ whose dimension is
length-squared. We take as our base shape $C$ a strip of length $L$, where
$L \gg \sqrt{A}$, and width $w$, where $w \ll \sqrt{A}$. We round the two
ends with semi-circles, as shown in on the top in figure 1.

As functions of a normal vector field $a.n$ along $C$, the metrics we want
to compare are:
\begin{enumerate}
\item $G^A_C(a,a) = \int_C (1+A\ka^2)a^2.ds$
\item $G^{\on{imm},1}_C(a,a) = \underset{b}{\inf} \int_C\Big(|a\cdot n +
b\cdot v)|^2 +A|D_s(a\cdot n + b\cdot v)|^2\Big) ds$
\item $G^{\on{diff},1}_C(a,a) = \tfrac1{\sqrt{A}}
\underset{\left(\substack{ { \mathbb R^2 \text{-vec.flds.}X}\\ { \langle
X,n\rangle = a }} \right)}{\inf} \iint_{\mathbb
R^2}\Big(|X|^2+A|DX|^2\Big)\,dx\,dy$,
\item $G^{\on{diff},2}_C(a,a) = \tfrac1{\sqrt{A}}
\underset{\left(\substack{{\mathbb R^2 \text{-vec.flds.}X}\\ { \langle
X,n\rangle = a}}\right)}{\inf} \iint_{\mathbb R^2}
\Big(|X|^2+2A|DX|^2+A^2|D^2 X|^2\Big)\,dx\,dy$
\end{enumerate}
The term $\tfrac1{\sqrt{A}}$ in the last 2 metrics is put there so that the
double integrals have the same `dimension' as the single integrals. In this
way, all the metrics will be comparable.

\begin{figure} \begin{center}
\psfrag{C}{$C$}
\psfrag{ph-}{$\ph_-$}
\psfrag{phf}{$\ph_f$}
\psfrag{ph+}{$\ph_+$}
\psfrag{phx}{$\ph_x$}
\psfrag{L}{$L$}
\psfrag{I}{$I$}
\psfrag{2w}{$2w$}
\psfrag{involution}{involution}
\epsfig{width=6cm,file=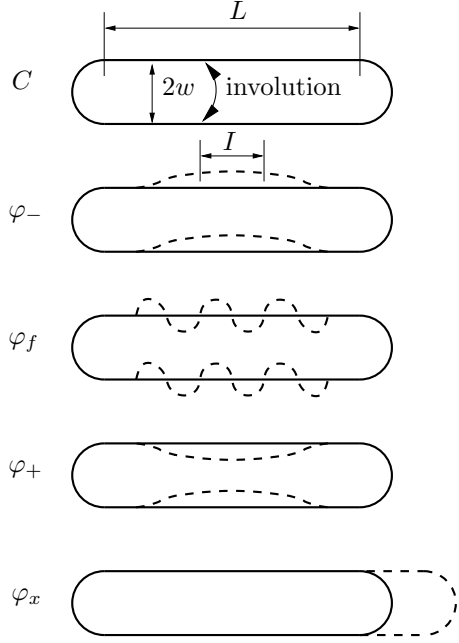}
\end{center}
\label{figure:cigar}
\caption{The cigarlike shape and their deformations}
\end{figure}
To compare the 4 metrics, we don't take all normal vector fields $a.n$
along $C$. Note that $C$ has an involution, which flips the top and bottom
edges. Thus we have even normal vector fields and odd normal vector fields.
Examples are shown in figure 1. We will consider two even and two odd
normal vector field, described below, and normalize each of them so that
$\int_C a^2ds = 1$. They are also shown in figure 1. 
They involve some
interval $I$ along the long axis of the shape of length $\la \gg w$. The
interval determines a part $I_t$ of the top part of $C$ and $I_b$ of the
bottom.
\begin{enumerate}
\item Let $a \equiv +1/\sqrt{2\la}$ along $I_t$ and $a \equiv
-1/\sqrt{2\la}$ along $I_b$, $a$ zero elsewhere except we smooth it at the
endpoints of $I$. Call this odd vector field $\ph_-$.
\item  Fix a high frequency $f$ and, on the same intervals, let $a(x) = \pm
\sin(f.x)/\sqrt{\la}$. Call this odd vector field $\ph_f$.
\item The third vector field is even and is defined by
$a(x)=\sqrt{\tfrac2{\pi w}}\langle n, \tfrac{\partial}{\partial x}\rangle$
at the right end of the curve, being zero along top, bottom and left end.
Call this $\ph_x$. The factor in front normalizes this vector field so that
its $L^2$ norm is 1.
\item Finally, we define another even vector field by $a=+1/\sqrt{2\la}$ on
both $I_t$ and $I_b$, zero elsewhere except for being smoothed at the ends
of $I$. Call this $\ph_+$.
\end{enumerate}

The following table shows the approximate leading term in the norm of each
of these normal vector fields $a$ in each of the metrics. By approximate,
we mean the exact norm is bounded above and below by the entry times a
constant depending only on $C$ and by leading term, we mean we ignore terms
in the small ratios $w/\sqrt{A}$, $\la/\sqrt{A}$:
$$ \begin{array}{ccccc}
\text{function} & G^A & G^{\on{imm},1} & G^{\on{diff},1} & G^{\on{diff},2}
\\
\hline
\ph_- & 1 & 1 & 1 & 1\\
\ph_f & 1 & (\sqrt{A} f)^2 & \sqrt{A}f & (\sqrt{A}f)^3\\
\ph_+ & 1 & 1 & \sqrt{A}/w & (\sqrt{A}/w)^3\\
\ph_x & A/w^2 & \sqrt{A}/w & \tfrac{\sqrt{A}}{w\log(A/w^2)} & \sqrt{A}/w
\end{array} $$
Thus, for instance:
\begin{align*}
 G^A(\ph_x,\ph_x) &= \frac2{\pi w} \int_{\text{right end}}(1+A\ka^2)
 \langle n, \tfrac{\partial}{\partial x}\rangle^2 ds \\
 &= \tfrac2{\pi w} (1+Aw^{-2})\ell(\text{right  end}) \text{Ave}(\langle n,
 \tfrac{\partial}{\partial x}\rangle^2)\\
 &= (1+Aw^{-2}) \approx A/w^2
 \end{align*}
The values of all the other entries under $G^A$ are clear because $\ka
\equiv 0$ in their support.

To estimate the other entries, we need to estimate the horizontal lift,
i.e., the functions $b$ or $v$.
To estimate the norms for $G^{\on{imm},1}$, we take $b=0$ in all cases
except $\ph_x$ and then get
$$ G^{\on{imm},1}_C(a,a) = \|a\|^2_{H^1_A}$$
the first Sobolev norm. For $a=\ph_f$, we simplify this, replacing the full
norm by the leading term $A(D_s(a))^2$ and working this out. To compute
$G^{\on{imm},1}(\ph_x,\ph_x)$, let $k=\sqrt{\tfrac2{\pi w}}$ be the
normalizing factor and lift $a.n$ along the right end of $C$ to the
$\mathbb R^2$ vector field $k.\tfrac{\partial}{\partial x}$. This adds a
tangential component which we taper to zero on the top and bottom of $C$
like $k.e^{-x/\sqrt{A}}$. This gives the estimate in the table.

Finally, consider the 2 metrics $G^{\on{diff},k}$, $k=1,2$. For these, we
need to lift the normal vector fields along $C$ to vector fields on all of
$\mathbb R^2$. For the two odd vector fields $f=\ph_-$ and $f=\ph_f$, we
take $v$ to be constant along the small vertical lines inside $C$ and zero
in the extended strip $-w \le y \le w, x \notin I$ and we define $v$
outside $-w \le y \le w$ by:
\begin{align*}
v(x,y+w) &= v(x,-y-w) = F(x,y) \tfrac{\partial}{\partial y}, \\
\widehat{F}(\xi,\et)&=\frac{k\sqrt{A}(1+A\xi^2)^{k-1/2}.\hat{f}(\xi)}{\pi(1
+A(\xi^2+\et^2))^k}
\end{align*}
We check the following:
\begin{align*}
&\text{(a)} \qquad \left((I-A\triangle)^k F\right)^\wedge =
\tfrac{k}{\pi}\sqrt{A}(1+A\xi^2)^{k-1/2} \hat{f}(\xi) \text{ is indep. of }
\et \text{ hence} \\
&\text{support}(I-A\triangle)^kF) \subset \{y=0\} \\
& \text{(b)} \qquad \int \hat{F} d\et = \hat{f}, \text{ hence } F|_{y=0} =
f.
\end{align*}
Thus:
\begin{align*}
G^{\on{diff},k}_C(f,f) &\approx \tfrac1{\sqrt{A}} \iint_{\mathbb R^2}
\langle (I-A\triangle)^kF,F \rangle dx dy \\
&=  \tfrac{k}{\pi} \int (1+A\xi^2)^{k-1/2} |\hat{f}(\xi)|^2 d\xi =
\tfrac{k}{\pi} \|f\|^2_{H^{k-1/2}_A}
\end{align*}
The leading term in the $k^{th}$ Sobolev norm of $\ph_f$ is
$(\sqrt{A}f)^{2k}$, which gives these entries in the table.

To estimate $G^{\on{diff},1}_C(\ph_+,\ph_+)$, we define $v$ by extending
$\ph_+$ linearly across the vertical line segments $-w \le y \le w, x \in
I$, i.e.\ to $\ph_+(x) y/w$. This gives the leading term now, as the
derivative there is $\ph_+/w$. In fact for any odd vector field $a$ of
$L^2$-norm 1 and for which the derivatives are sufficiently small compared
to $w$, the norm has the same leading term:
$$G^{\on{diff},1}_C(a,a) \approx \sqrt{A}/w.$$

To estimate $G^{\on{diff},2}_C(\ph_+,\ph_+)$, we need a smoother extension
across the interior of $C$. We can take $\ph_+(x).\tfrac32(\tfrac{y}{w} -
\tfrac13 (\tfrac{y}{w})^3).$ Computing the square integral of the second
derivative, we get the table entry $G^{\on{diff},2}_C(a,a) \approx
(\sqrt{A}/w)^3.$

To estimate $G^{\on{diff},k+1}_C(\ph_x,\ph_x)$, we now take $v$ to be
$$v=c(k,A,w) \left[ \left(\tfrac{|x|}{\sqrt{A}}\right)^k
K_k\left(\tfrac{|x|}{\sqrt{A}}\right) \ast \ch_D \right]
\frac{\partial}{\partial x}.$$
where $D$ is the disk of radius $w$ containing the arc making up the right
hand end of $C$, and where $c(k,A,w)$ is a constant to be specified later. 
The function 
$\tfrac1{2\pi k! A} \left(\tfrac{|x|}{\sqrt{A}}\right)^k
K_k\left(\tfrac{|x|}{\sqrt{A}}\right)$ 
is the fundamental solution of
$(I-A\triangle)^{k+1}$ and is $C^1$ for $k>0$ but with a log pole at 0 for
$k=0$. Thus:
$$ (I-A\triangle)^{k+1}v = 2\pi k! Ac(k,A,w)\chi_D.\frac{\partial}{\partial
x}$$
while, up to upper and lower bounds depending only on $k$, the restriction
of $v$ to the disk $D$ itself is equal to $\log(\sqrt{A}/w)c(0,A,w)w^2$ if
$k=0$ and simply $c(k,A,w)w^2$ for $k>0$. By symmetry $v$ is also constant
on the boundary of $D$ and thus $v$ extends $\ph_x$ if we take $c(0,A,w) =
c_0/\log(\sqrt{A}/w)w^{5/2}$ if $k=0$ and $c(k,A,w)=c_k/w^{5/2}$ if $k>0$
(constants $c_k$ depending only on $k$). Computing the $H^k$-norm of $v$,
we get the last table entries.

Summarizing, we can say that the large norm of $\ph_x$ is what
characterizes $G^A$; the large norms of $\ph_+$ characterize
$G^{\on{diff}}$; and the rate of growth in frequency of the norm of $\ph_f$
distinguishes all 4 norms.

\bibliographystyle{plain}

\end{document}